\documentclass[11pt]{article}
\usepackage[english]{babel}
\usepackage{amsmath,amssymb}
\usepackage{mathtools}
\usepackage{bbm}
\usepackage{graphicx}
\usepackage{color}
\usepackage{multirow}
\usepackage{tikz}
\numberwithin{equation}{section}
\newcounter{todocounter}

\usepackage{float}

\usepackage{xargs} 

\usepackage{accents}

\usepackage{enumerate}

\usepackage{tabulary,capt-of}

\usepackage{ctable}

\usepackage{subfig}

\usepackage{url}

\usepackage[hidelinks]{hyperref}

\usepackage[toc,page]{appendix}

\usepackage[colorinlistoftodos,prependcaption,textsize=tiny]{todonotes}
\newcommandx{\unsure}[2][1=]{\todo[linecolor=red,backgroundcolor=red!25,bordercolor=red,#1]{#2}}
\newcommandx{\change}[2][1=]{\todo[linecolor=blue,backgroundcolor=blue!25,bordercolor=blue,#1]{#2}}
\newcommandx{\info}[2][1=]{\todo[linecolor=OliveGreen,backgroundcolor=OliveGreen!25,bordercolor=OliveGreen,#1]{#2}}
\newcommandx{\improvement}[2][1=]{\todo[linecolor=yellow,backgroundcolor=yellow!25,bordercolor=yellow,#1]{#2}}
\newcommandx{\thiswillnotshow}[2][1=]{\todo[disable,#1]{#2}}
	
\usepackage[ulem=normalem]{changes}
\definecolor{darkgreen}{rgb}{0.0, 0.4, 0.0}
\definechangesauthor[name={Claudius Birke},color=red]{AT}
\setauthormarkup{}    
	
\title{
	\bf{A low Mach two-speed relaxation scheme for the compressible Euler equations with gravity}
}
\author{
	Claudius Birke\footnote{Fakult\"at f\"ur Mathematik und Informatik, Universit\"at W\"urzburg, Emil-Fischer-Str. 40, 97074 W\"urzburg, Germany} \footnote{Correspondence to: Claudius Birke, Universit\"at W\"urzburg, Emil-Fischer-Str. 40, 97074 W\"urzburg, Germany, Email: claudius.birke@mathematik.uni-wuerzburg.de},
Christophe Chalons\footnote{Laboratoire de Math\'{e}matiques de Versailles, UVSQ, CNRS, Universit\'{e} Paris-Saclay, 78035 Versailles, France},\\ Christian Klingenberg\footnotemark[1]
}
\date{\today}
	
\DeclareCaptionLabelSeparator{none}{ }

\begin{document}

\maketitle

\section*{Abstract}
We present a numerical approximation of the solutions of the Euler equations with a gravitational source term. On the basis of a Suliciu type relaxation model with two relaxation speeds, we construct an approximate Riemann solver, which is used in a first order Godunov-type finite volume scheme.
This scheme can preserve both stationary solutions and the low Mach limit to the corresponding incompressible equations. In addition, we prove that our scheme preserves the positivity of density and internal energy, that it is entropy satisfying and also guarantees not to give rise to numerical checkerboard modes in the incompressible limit. 
Later we give an extension to second order that preserves positivity, asymptotic-preserving and well-balancing properties. 
Finally, the theoretical properties are investigated in numerical experiments.

\paragraph{Keywords} Euler equations, finite volume methods, relaxation, well-balancing, low Mach, asymptotic-preserving, entropy satisfying, checkerboard modes, positivity preserving

\paragraph{AMS subject classification} 65M08, 76M12

\newtheorem{theorem}{Theorem}
\newtheorem{example}{Example}
\newtheorem{lemma}{Lemma}
\newtheorem{remark}{Remark}
\newtheorem{assumption}{Assumption}

\newcommand\mytodos[1]{\textcolor{red}{#1}}
\captionsetup[figure]{labelfont={bf},labelformat={default},labelsep=period,name={Fig.},labelsep=none}
\captionsetup[table]{labelfont={bf},labelsep=none}


\section{Introduction}
\label{sec:introduction}

The goal of this paper is to find a numerical approximation of the solutions of the Euler equations including a gravitational source term. In a dimensionless form these equations are defined by
\begin{equation}
\label{sys:Euler}
\begin{split}
\partial_t \rho + \nabla \cdot \left(\rho \boldsymbol{u} \right) &= 0, \\
\partial_t \left( \rho \boldsymbol{u} \right) + \nabla \cdot \left(\rho \boldsymbol{u} \otimes \boldsymbol{u} \right) + \frac{1}{M^2 }\nabla p  &= -\frac{1}{M^2}\rho \nabla \Phi, \\
\partial_t  E + \nabla \cdot \left( (E+p) \boldsymbol{u} \right) &= -\rho \boldsymbol{u} \cdot \nabla \Phi,
\end{split}
\end{equation}
where $\rho(\boldsymbol{x},t):\mathbb{R}^d \times \mathbb{R}_{\geq 0} \rightarrow \mathbb{R}^+$ denotes the density, $\boldsymbol{u}(\boldsymbol{x},t):\mathbb{R}^d \times \mathbb{R}_{\geq 0} \rightarrow \mathbb{R}^d$ the velocity vector, $E(\boldsymbol{x},t):\mathbb{R}^d \times \mathbb{R}_{\geq 0} \rightarrow \mathbb{R}^+$ the total energy and $\Phi(\boldsymbol{x}): \mathbb{R}^d \rightarrow \mathbb{R}$ a given smooth gravitational potential. In this dimensionless formulation, the parameter $M$ represents the Mach number, which controls the ratio between the velocity of the gas and the speed of sound.
In this work, we consider the combined low Mach/low Froude number limit, which is the reason why we set $Fr = M$. As an effect, only the Mach number $M$ appears in the dimensionless equations in (\ref{sys:Euler}). 
\\
The pressure is given by a pressure law $p(\tau,e): \mathbb{R}^+ \times \mathbb{R}^+ \rightarrow \mathbb{R}$, where $\tau=1/\rho$ denotes the specific volume and $e>0$ the internal energy. The total energy can then be expressed by 
\begin{equation}
\label{eq:energy}
E = \rho e + \frac{1}{2} M^2 \rho \left \lvert \boldsymbol{u} \right \rvert^2.
\end{equation} 
The pressure law closing this model obeys the second law of thermodynamics so that a specific entropy $s(\tau,e): \mathbb{R}^+ \times \mathbb{R}^+ \rightarrow \mathbb{R^+}$, which satisfies the relation
\begin{equation}
\label{eq:temperature}
-T ds = de + pd\tau
\end{equation}
for some temperature $T(\tau,e) > 0$, exists. In this work, we assume $(\tau,e) \mapsto s(\tau,e)$ to be strictly convex. \\
The phase space to which the system (\ref{sys:Euler}) is associated is denoted by
\begin{equation}
\label{eq:phase_space}
\Omega = \{ \omega = (\rho, \rho \boldsymbol{u}, E)^T \in \mathbb{R}^{d+2}; \rho >0, e>0  \}.
\end{equation} \\
This model can be used in various fields of application, such as the simulation of gas flows in the interior of stars in astrophysics. Depending on the application, the flows can have large scale differences, e.g. the sound speed can be much higher than the speed of the fluid flow.
In these low Mach number regimes standard finite volume schemes suffer from excessive diffusion, which 
can erase the structure of the solution beyond recognition.
A number of different strategies have been developed to overcome this problem.
One simple but efficient strategy is to modify the diffusion term in the numerical flux by rescaling it with the local Mach number and thereby reduce the viscosity on the velocity.
First introduced for the homogeneous Euler equations \cite{6,9,10,14,15,16,17,18}, this approach was also extended to the Euler equations including a gravitational source term \cite{19}.
A second approach introduced by Klein \cite{11} relies on a pressure splitting, which decomposes the system of equations into one slow, non-linear part and into a linear part for the fast acoustic dynamics. 
In \cite{34} this splitting is combined with a Suliciu type relaxation model and an implicit time integration. Thomann et al \textit{et al.} modify the approach to an implicit-explicit (IMEX) scheme, in which only the acoustic part is solved implicitly, while the non-linear part is solved explicitly by a Godunov-type method based on an approximate Riemann solver. Later this IMEX approach was extended to the Euler equations with gravity \cite{8}.

Basis of the herein presented scheme is a third alternative introduced by Chalons \textit{et al.} in \cite{1}, where a two-speed relaxation scheme for the barotropic, homogeneous Euler equations is proposed. The use of two different relaxation speeds enables an independent control of the numerical viscosity on the density and on the velocity. By special definitions of the speeds in the low Mach regime, viscosity is transferred from the velocity to the density. Therefore, this approach is related to the previously described rescaling of the viscosity term. The key advantage of this method is that under a subcharacteristic condition it is stable and provably entropy satisfying.
Later the two-speed relaxation system was used to develop an IMEX scheme for the homogeneous Euler equations \cite{3}.\\ 
In contrast, the method presented in this paper is fully explicit. The basic structure of the two-speed relaxation system is adopted and extended by gravitational source terms. 
From the exact resolution of the Riemann problem associated with this relaxation system, a Godunov-type finite volume method is constructed. The modification of the relaxation speeds is adopted from the original approach. The resulting approximate Riemann solver satisfies a discrete entropy inequality. 
Based on this inequality, it is shown that no checkerboard modes can arise in the variables fluid velocity and pressure.
Checkerboard modes pose an instability characterized by a decoupling of the spatial approximation, which can occur in numerical solutions of incompressible fluid equations computed on collocated grids \cite{31}.
It is well-known that most of the asymptotic-preserving schemes exhibit such nonphysical checkerboard modes in low Mach regimes \cite{13,17}.

When studying the Euler equations with gravity source terms one has to consider their influence on the behaviour of steady states.
In several applications such as astrophysics one deals with problems close to the hydrostatic equilibrium
\begin{align}
\label{sys:hydrostatic_equilibrium}
\left\{\begin{array}{ll} \boldsymbol{u}=0, \\
							  \nabla p = -\rho \nabla \Phi.
										\end{array}\right.
\end{align}
Standard finite volume schemes do not automatically satisfy a discrete equivalent of (\ref{sys:hydrostatic_equilibrium}). Therefore these steady states are not preserved exactly by such schemes and small perturbations around this equilibrium cannot be resolved unless the resolution of the scheme is increased, so that the truncation error is sufficiently small. In order to avoid this potentially high computational effort, well-balanced schemes \cite{4,5,8,20,21,27,28,29,32} were introduced, which satisfy exactly a discrete equivalent of the steady state.

The well-balancing mechanism in the herein presented relaxation scheme is taken over from \cite{2}. The key idea is to add a transport relaxation equation for the gravitational potential to the relaxation system, which leads
to a Riemann-problem that  is under-determined. This gives an additional degree of freedom and allows to introduce a closure equation that is a discrete equivalent of (\ref{sys:hydrostatic_equilibrium}) and ensures the well-balanced property. 
This approach is exact for certain families of hydrostatic equilibria, i.e. isothermal, incompressible and polytropic ones. In all other cases it maintains the equilibrium to second order. We extend this approach so that it can be applied to any hydrostatic solution for the Euler equations with any equation of state if the hydrostatic solution is known a priori. The extension is based on a second order approximation of the difference in the gravitational potential using the given hydrostatic states for the density and pressure. This is useful for applications in stellar astrophysics, in which the equation of state (EoS) is given in form of a table. Since hydrostatic solutions depend on the EoS, they can then only be found through numerical simulations carried out beforehand and are therefore available in the form of discrete data.

The paper is organized as follows. In Sect. \ref{sec:relaxation_model}, the two-speed relaxation model is derived. In addition, the approximate Riemann solver associated with this system and its intermediate states are determined. The following Sect. \ref{sec:relaxation_scheme} contains the first order Godunov-type finite volume scheme, which is based on the previously introduced approximate Riemann solver. Its properties are described and proven in Sect. \ref{sec:properties}. A suitable extension to second order in space is given in Sect. \ref{sec:2nd_order}. In Sect. \ref{sec:results}, the properties of the second order scheme are checked in numerical tests. 
Finally, Sect. \ref{sec:conclusion} provides the conclusion and an outlook.


\section{The relaxation model}
\label{sec:relaxation_model}

The one dimensional relaxation system described below is based at its core on the Suliciu relaxation model \cite{22,23,24}. The pressure $p$ is approximated by the relaxation variable $\pi$ and we add an additional equation describing its behaviour to the system
\begin{equation}
\label{eq:relaxation_pi_equation}
\partial_t \rho \pi + \partial_x (\rho \pi v) + ab \partial_x v = \rho \frac{p-\pi}{\varepsilon}.
\end{equation}
While only one relaxation speed is used in the classical Suliciu relaxation model, here two speeds $a>0$ and $b>0$ appear, as proposed in \cite{1}. This will be useful to control viscosity for pressure and velocity separately. These speeds will be defined later in Sect. \ref{sec:AP} so that they meet stability criteria and keep the viscosity bounded in the low Mach regime.

In addition, also the velocity $u$ is approximated by a relaxation variable $v$ and the following equation is introduced
\begin{equation}
\label{eq:relaxation_v_equation}
\partial_t \left( \rho v \right) + \partial_x \left(\rho v^2 \right) + \frac{a}{b} \partial_x \frac{\pi}{M^2} =  \rho \frac{u-v}{\varepsilon} -\frac{a}{b} \frac{1}{M^2} \rho \partial_x \Phi.
\end{equation}
In the next step, we also want to include the gravitational potential in the approximate Riemann solver. According to \cite{2} this can be done by approximating the gravitational potential $\Phi$ by the relaxation variable $Z$ and adding a transport relaxation equation to the relaxation system
\begin{equation}
\label{eq:relaxation_Z_equation}
\partial_t \rho Z + \partial_x \rho v Z = \rho \frac{\Phi-Z}{\varepsilon}.
\end{equation}
Finally, we derive the following relaxation model
\begin{equation}
\label{sys:relaxation_model}
\begin{split}
\partial_t \rho + \partial_x \left(\rho v \right) &= 0, \\
\partial_t \left( \rho u \right) + \partial_x \left(\rho uv+\frac{\pi}{M^2} \right) &= - \frac{1}{M^2} \rho \partial_x Z, \\
\partial_t E + \partial_x \left( (E+\pi)v \right) &= - \rho v \partial_x Z, \\
\partial_t \left(\rho \pi \right) + \partial_x \left(\rho \pi v \right) + ab\partial_x v &= \rho \frac{p-\pi}{\varepsilon}, \\
\partial_t \left( \rho v \right) + \partial_x \left(\rho v^2 \right) + \frac{a}{b} \partial_x \frac{\pi}{M^2} &=  \rho \frac{u-v}{\varepsilon} -\frac{a}{b} \frac{1}{M^2} \rho \partial_x Z, \\
\partial_t \rho Z + \partial_x \rho v Z &= \rho \frac{\Phi-Z}{\varepsilon}, \\
\partial_t a + v\partial_x a &= 0, \\
\partial_t b + v\partial_x b &= 0.
\end{split}
\end{equation}
The solutions to this relaxation model can be seen as a viscous approximation of the solutions of the original system (\ref{sys:Euler}) as long as the subcharacteristic conditions 
\begin{equation}
\label{eq:subcharacteristic_condition}
a \geq b \quad \text{and} \quad ab \geq \rho^2 c^2
\end{equation}
are satisfied.\\

\begin{remark}
\label{rem:Suliciu_model}
By choosing $u=v$ and $a=b$ one recovers the standard Suliciu relaxation model.
\end{remark}

The homogeneous system, denoted by (\ref{sys:relaxation_model})$_{\varepsilon=\infty}$, has the following properties.\\

\setcounter{lemma}{1}
\begin{lemma}
\label{lem:wave_speeds_and_riemann_invariants}
The relaxation system (\ref{sys:relaxation_model})$_{\varepsilon=\infty}$ is hyperbolic and all characteristic fields are linearly degenerate. The eigenvalues of the system are given by
\begin{equation}
\label{eq:wave_speeds}
\sigma^v = v, \quad \sigma^{\pm}=v \pm \frac{a}{M \rho}
\end{equation}
where $\sigma^v$ has multiplicity six. The eigenvalues have the fixed ordering
\begin{equation}
\label{eq:order_wave_speeds}
\sigma^-<\sigma^v<\sigma^+.
\end{equation}
The Riemann invariant corresponding to the eigenvalue $\sigma^v$ is
\begin{equation}
\label{eq:riemann_invariant_v}
I_1^{v}=v
\end{equation}
and those corresponding to $\sigma^{\pm}$ are
\begin{equation}
\label{eq:riemann_invariant_+-}
\begin{split}
I_1^{\pm} &= v \pm \frac{a}{M\rho}, \  \ I_2^{\pm}=u \pm \frac{b}{M\rho}, \ \  I_3^{\pm}=  \frac{1}{\rho}+\frac{\pi}{ab}, \\
I_4^{\pm} &= e+\frac{(a-b)b+2\rho(\pi\pm bM(v-u))}{2\rho^2}, \\
 I_5^{\pm} &= a, \ \ I_6^{\pm}=b, \ \  I_7^{\pm}=Z.
\end{split}
\end{equation}
\end{lemma}

\begin{proof}
The computations are straightforward and left to the reader.
\hfill
\end{proof} 

\setcounter{remark}{2}
\begin{remark}
\label{rem:under-determined_riemann_problem}
The relaxation system (\ref{sys:relaxation_model})$_{\varepsilon=\infty}$ provides only one Riemann invariant $I_1^v$ for the contact wave. As a result, the associated Riemann problem is under-determined.
\end{remark} 

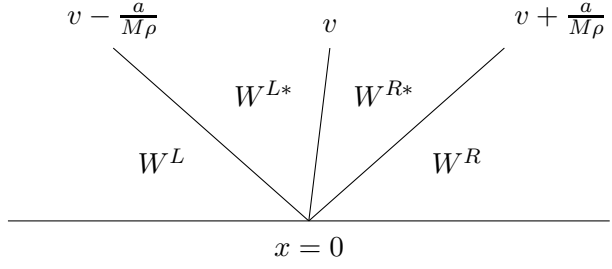
\begin{figure}[t]
	\centering
	\begin{tikzpicture}[scale=2.]
	\begin{scope}[every node/.style={scale=1.}]
	\draw
	(-1.0,0)-- (3.0,0);
	\draw (-0.3,1.15) coordinate (l1) node[above] {$v - \frac{a}{M\rho}$} -- (1,0) coordinate (x0) node[below, yshift=-0.5ex] {$x=0$} ;
	\draw (2.3,1.15) coordinate (l3) node[above right] {$v + \frac{a}{M\rho}$} -- (1,0);
	\draw (1.14,1.15) coordinate (l2) node[above,yshift=0.5ex] {$v$} -- (1,0);
	\draw (0.25,0.4) coordinate (WL) node[left] {$W^L$};
	\draw (1.75,0.4) coordinate (WR) node[right] {$W^R$};
	\draw (0.7,0.7) coordinate (WSL) node[above] {$W^{L*}$};
	\draw (1.5,0.7) coordinate (WSR) node[above] {$W^{R*}$};
	\end{scope}
	\end{tikzpicture}
	\caption{Schematic diagram of the Riemann fan for the relaxation system (\ref{sys:relaxation_model}). 
The Riemann fan consists of the
two intermediate states $W^{L*}$ and $W^{R*}$ for given states $W^L$ and $W^R$. The states are separated by the three wave speeds $v-a/(M\rho)$, $v$ and $v+a/(M\rho)$}
\label{fig:riemann_fan}
\end{figure}

Let us now consider a single Riemann problem associated with the system (\ref{sys:relaxation_model})$_{\varepsilon=\infty}$. In order to simplify the notations we introduce the state vector
\begin{equation}
\label{eq:model_state_vector}
W = \left( \rho, \rho u, E, \rho \pi, \rho v, \rho Z,a,b \right)^T
\end{equation} 
in the phase space
\begin{equation}
\mathcal{O} = \{ W\in \mathbb{R}^8: \ \rho > 0, \ e>0 \}
\end{equation}
and additionally for $\omega \in \Omega$ and given gravitational potential $\Phi$ the state vector at relaxation equilibrium denoted by
\begin{equation}
W^{eq}(\omega) =  \left( \rho, \rho u, E, \rho p(\tau,e), \rho u, \rho \Phi,a,b \right)^T.
\end{equation}
Then the initial data of the Riemann problem is given by two constant states $W^L$ and $W^R$ separated by one discontinuity located at $x=0$
\begin{equation}
\label{eq:riemann_problem}
W_0(x) = 
\left\{\begin{array}{ll} W^L, & x<0, \\
							  W^R, & x>0.
										\end{array}\right.
\end{equation}
The solution to this problem consists of four constant states, each separated by a contact discontinuity. Therefore the approximate Riemann solver $W_{\mathcal{R}}(x/t;W^L,W^R)$ has the structure
\begin{equation}
\label{eq:approximate_riemann_solver}
W_{\mathcal{R}}(\frac{x}{t};W^L,W^R) = 
\left\{\begin{array}{ll} W^L, & \frac{x}{t}<\sigma^-, \\
							  W^{L*}, & \sigma^-<\frac{x}{t}<\sigma^v, \\
							  W^{R*}, & \sigma^v<\frac{x}{t}<\sigma^+, \\
							  W^R, & \sigma^+<\frac{x}{t}.
										\end{array}\right.
\end{equation}
This structure of the solution is also shown in Fig. \ref{fig:riemann_fan}.
For the computation of the intermediate states $W^{L*}$ and $W^{R*}$ we can use the Riemann invariants given in lemma \ref{lem:wave_speeds_and_riemann_invariants}. Since Riemann invariants are constant across their corresponding wave, each Riemann invariant provides one equation. However, counting the Riemann invariants reveals that only 15 Riemann invariants face 16 unknown intermediate states. Therefore, the Riemann problem (\ref{eq:riemann_problem}) is, as already stated in remark \ref{rem:under-determined_riemann_problem}, under-determined. In order to overcome this problem, it is suggested in \cite{2} to introduce an additional relation
\begin{equation}
\label{eq:closure_equation}
\pi^{R*}-\pi^{L*} = -\bar{\rho} \left(W^L,W^R \right) \left( Z^R-Z^L \right),
\end{equation}
where the function $\bar{\rho}$ denotes a $\rho$-average function. This equation is chosen because it is a discrete representation of the steady states at rest in (\ref{sys:hydrostatic_equilibrium}) in one spatial dimension and therefore will be useful for the well-balancing of hydrostatic equilibria. The explicit definition of the function $\bar{\rho}$ depends on the underlying hydrostatic equilibrium and will be given later in Sect. \ref{sec:wb}.

With the newly added closure equation, it is now possible to compute the intermediate states in the Riemann solution.\\

\setcounter{lemma}{3}
\begin{lemma}
The solution of the Riemann problem (\ref{eq:riemann_problem}) associated with the relaxation system (\ref{sys:relaxation_model})$_{\varepsilon=\infty}$ has the structure given in (\ref{eq:approximate_riemann_solver}) with the intermediate states
\begingroup
\allowdisplaybreaks
\begin{align}
\label{eq:IS_vs}
&v^* = \frac{M b^L v^L + M b^R v^R + \pi^L-\pi^R - \bar{\rho} \left(W^L,W^R \right) \left( Z^R-Z^L \right)}{M(b^L+b^R)}, \\
\label{eq:IS_rhoLs}
&\frac{1}{\rho^{L*}} = \frac{1}{\rho^L} + \frac{M b^R \left( v^R-v^L \right) + \pi^L-\pi^R - \bar{\rho} \left(W^L,W^R \right) \left( Z^R-Z^L \right)}{a^L \left( b^L+b^R \right)}, \\
\label{eq:IS_rhoRs}
&\frac{1}{\rho^{R*}} = \frac{1}{\rho^R} + \frac{M b^L \left( v^R-v^L \right) + \pi^R-\pi^L + \bar{\rho} \left(W^L,W^R \right) \left( Z^R-Z^L \right)}{a^R \left( b^L+b^R \right)}, \\
\label{eq:IS_uLs}
&u^{L*} = u^L + \frac{b^L \left( b^R M \left(v^R-v^L \right) + \pi^L-\pi^R - \bar{\rho} \left(W^L,W^R \right) \left( Z^R-Z^L \right) \right)}{M a^L \left( b^L+b^R \right) }, \\
\label{eq:IS_uRs}
&u^{R*} = u^R+ \frac{b^R \left( b^L M \left(v^L-v^R \right) + \pi^L-\pi^R - \bar{\rho} \left(W^L,W^R \right) \left( Z^R-Z^L \right) \right)}{M a^R \left( b^L+b^R \right) }, \\
\label{eq:IS_piLs}
&\pi^{L*} = \frac{b^R \pi^L + b^L \pi^R + M b^L b^R \left(v^L-v^R \right) + b^L \bar{\rho} \left(W^L,W^R \right) \left( Z^R-Z^L \right) }{b^L+b^R}, \\
\label{eq:IS_piRs}
&\pi^{R*} = \frac{b^R \pi^L + b^L \pi^R + M b^L b^R \left(v^L-v^R \right) - b^R \bar{\rho} \left(W^L,W^R \right) \left( Z^R-Z^L \right) }{b^L+b^R}, \\ 
\label{eq:IS_eLs} 
&e^{L*} = e^L + \frac{(\pi^{L*})^2-(\pi^L)^2}{2 a^L b^L} + \frac{(v^*-u^{L*})^2-(v^L-u^L)^2}{2 (\frac{a^L}{b^L}-1)}, \\
\label{eq:IS_eRs}
&e^{R*} = e^R + \frac{(\pi^{R*})^2-(\pi^R)^2}{2 a^R b^R} + \frac{(v^*-u^{R*})^2-(v^R-u^R)^2}{2 (\frac{a^R}{b^R}-1)}, \\
\label{eq:IS_abZ}
&a^{L*} = a^L, \ a^{R*} = a^R, \ b^{L*} = b^L, \ b^{R*} = b^R, \ Z^{L*} = Z^L, \ Z^{R*} = Z^R.
\end{align}
\endgroup
\end{lemma}
\begin{proof}
The intermediate states can be computed by solving the system of equations given by the Riemann invariants and the closure equation (\ref{eq:closure_equation}). The precise steps are straightforward and therefore left to the reader.
\hfill
\end{proof}

\setcounter{remark}{4}
\begin{remark}
\label{rem:defintion_relaxation_speeds}
At this point, we do not explicitly define the relaxation speeds $a^L$, $a^R$, $b^L$ and $b^R$, since later, in the proofs of the properties of the relaxation method, various conditions are placed on these speeds. The explicit definitions are then provided in Sect. \ref{sec:AP}.
\end{remark}

Equipped with the approximate Riemann solver, we can now define the overall discretization of the scheme in the next section.


\section{The relaxation scheme}
\label{sec:relaxation_scheme}
Before we derive a complete finite volume scheme for the Euler equations with a gravitational source (\ref{sys:Euler}), we introduce some useful notations. The spatial domain is divided into cells $\mathcal{C}_i = (x_{i-1/2},x_{i+1/2})$ with $i\in \mathbb{Z}$ that have the size $\Delta x = x_{i+1/2}-x_{i-1/2}$. The cell centers are denoted by $x_i$. The time discretization is given by $t^n=n \Delta t$ with $n\in \mathbb{N}$ and a timestep $\Delta t$ that is restricted by the CFL condition
\begin{equation}
\label{eq:cfl_condition}
\frac{\Delta t}{\Delta x} \underaccent{i}{\max} \left \{  \left \lvert v_i - \frac{a_{i}}{M \rho_i} \right \rvert, \left \lvert v_i + \frac{a_{i}}{M \rho_i}  \right \rvert \right \} \leq \frac{1}{2}. 
\end{equation}
The cell average $\omega_i^n$ then approximates the value over the cell $\mathcal{C}_i$ at time $t^n$
\begin{equation}
\label{eq:cell_average}
\omega_i^n \approx \frac{1}{\Delta x} \int_{\mathcal{C}_i} \omega (x,t^n) dx.
\end{equation}

At the start of each time step, we assume to be at the relaxation equilibrium. Therefore the initial data for the relaxation variables at time level $n$ is defined by
\begin{equation}
\pi_i^n =  p_i^n, \quad v_i^n=u_i^n, \quad Z_i^n = \Phi_i^n.
\end{equation} 
Starting from the equilibrium we solve the homogeneous relaxation system (\ref{sys:relaxation_model})$_{\varepsilon=\infty}$ using the Riemann solver $W_{\mathcal{R}}$ defined in (\ref{eq:approximate_riemann_solver}) and update the cell averages to the next time level $t^{n+1}$ by a Godunov method of the form
\begin{equation}
\label{eq:relaxatation_scheme}
\begin{split}
\omega_i^{n+1} =& \omega_i^n - \frac{\Delta t}{\Delta x} \left( F_{i+1/2}^n - F_{i-1/2}^n \right)\\ 
&+ \frac{\Delta t}{2} 
\left( S_{i-1/2}^{+,n} \frac{\Phi_{i}^n-\Phi_{i-1}^n}{\Delta x} + S_{i+1/2}^{-,n} \frac{\Phi_{i+1}^n-\Phi_{i}^n}{\Delta x} \right), \\
F_{i-1/2}^n =& F (\omega_{i-1}^n,\Phi_{i-1}^n,\omega_{i}^n, \Phi_{i}^n), \quad \ \  F_{i+1/2}^n = F (\omega_i^n,\Phi_i^n,\omega_{i+1}^n, \Phi_{i+1}^n), \\ 
S_{i-1/2}^{+,n} =& S^+ (\omega_{i-1}^n,\Phi_{i-1}^n,\omega_{i}^n, \Phi_{i}^n), \quad 
S_{i+1/2}^{-,n} =   S^- (\omega_i^n,\Phi_i^n,\omega_{i+1}^n, \Phi_{i+1}^n).
\end{split}
\end{equation}
The numerical flux is defined by
\begin{equation}
\label{eq:numerical_flux}
F(\omega^L,\Phi^L,\omega^R,\Phi^R) = 
\left\{\begin{array}{ll} F(\omega^L), & \text{ if } \sigma^- > 0, \\
							  F^{L*}, & \text{ if } \sigma^- < 0 \leq \sigma^v, \\
							  F^{R*}, & \text{ if } \sigma^v < 0 < \sigma^{+}, \\
							  F(\omega^R), & \text{ if } \sigma^+ < 0,
										\end{array}\right.
\end{equation}
where according to the left-hand sides of the first three equations of (\ref{sys:relaxation_model}) the intermediate fluxes can be written as
\begin{equation}
\label{eq:intermediate_fluxes}
\begin{split}
F^{L*} &= \left( \rho^{L*} v^*, \rho^{L*} u^{L*} v^* + \frac{\pi^{L*}}{M^2}, (E^{L*}+\pi^{L*})v^* \right), \\
F^{R*} &= \left( \rho^{R*} v^*, \rho^{R*} u^{R*} v^* + \frac{\pi^{R*}}{M^2},(E^{R*}+\pi^{R*})v^* \right).
\end{split}
\end{equation}
The numerical source terms are set as follows
\begin{equation}
\label{eq:source_terms}
\begin{split}
S^+(\omega^L,\Phi^L,\omega^R,\Phi^R) &= -(\text{sgn} (v^*)+1) \left( 0, \frac{1}{M^2 }\bar{\rho} (W^L,W^R) , \bar{\rho} (W^L,W^R) v^* \right)^T, \\
S^-(\omega^L,\Phi^L,\omega^R,\Phi^R) &= \ \  (\text{sgn} (v^*)-1) \left( 0, \frac{1}{M^2} \bar{\rho} (W^L,W^R) , \bar{\rho} (W^L,W^R) v^* \right)^T.
\end{split}
\end{equation}
We note that in this procedure only the variables of the original Euler equations (\ref{sys:Euler}) in the vector $\omega$ are updated to the next time level. For the upcoming time step we again assume to be at the equilibrium. 
As a consequence of this projection approach, the relaxation parameter $\varepsilon$ does not appear in the relaxation scheme (\ref{eq:relaxatation_scheme}) and thus does not have to be set explicitly.


\section{Properties of the relaxation scheme}
\label{sec:properties}

In this section we focus on the properties of the relaxation scheme just described. We start with the property of entropy stability.

\subsection{Entropy inequality}
\label{sec:entropy}
We seek those correct solutions that satisfy the entropy inequality.
In practice, it can be observed that searching for entropy solutions makes a finite volume method more stable. This is partly because an entropy inequality can help to ensure the positivity of density and/or internal energy.

Going back to the Euler equations (\ref{sys:Euler}) and assuming smooth solutions, it is possible to derive the additional conservation law
\begin{equation}
\label{eq:entropy_conservation}
\partial_t \rho \mathcal{F}(s) + \partial_x \rho \mathcal{F}(s) u = 0
\end{equation}
for all smooth functions $\mathcal{F}$. Assuming that $\mathcal{F}$ is increasing and $\omega \mapsto \rho \mathcal{F}(s)$ is convex, the pair $(\rho \mathcal{F}(s),\rho \mathcal{F}(s)u)$ defines a Lax entropy pair for the system \eqref{sys:Euler}.
Thus, equation \eqref{eq:entropy_conservation} states that the entropy is conserved for smooth solutions. However, since the Euler equations are non-linear, discontinuities can arise in the solution in finite time despite of smooth initial conditions. 
At discontinuities the equation (\ref{eq:entropy_conservation}) is not valid, since it does not consider the entropy dissipation at shocks. Therefore, we replace the equality in (\ref{eq:entropy_conservation}) by an inequality, which leads to the following entropy inequality
\begin{equation}
\label{eq:continuous_entropy_inequality}
\partial_t \rho \mathcal{F}(s) + \partial_x \rho \mathcal{F}(s) u \leq 0.
\end{equation}
Our scheme should now mimic this behaviour in the sense that its solutions satisfy a discrete version of (\ref{eq:continuous_entropy_inequality}).

\begin{theorem}
\label{theo:entropy}
Let us assume that $w_i^n$ belongs to $\Omega$ for all $i \in \mathbb{Z}$.
Furthermore, we assume that at each interface with initial left state $\omega^L$ and initial right state $\omega^R$  
the intermediate states for density and internal energy in the Riemann solution are positive, i.e. $\rho^{L*},\rho^{R*},e^{L*},e^{R*}>0$, and that the relaxation speeds $a^{L,R}$ and $b^{L,R}$ are such that they satisfy the subcharacteristic Whitham conditions
\begin{align}
\label{eq:whitham_condition1}
a^L b^L &> p(\tau^L,e^L) \partial_e p(\tau^L,e^L) - \partial_\tau p(\tau^L,e^L), \\
\label{eq:whitham_condition2}
a^L b^L &> p(\tau^{L*},e^{L*}) \partial_e p(\tau^{L*},e^{L*}) - \partial_\tau p(\tau^{L*},e^{L*}), \\
\label{eq:whitham_condition3}
a^R b^R &> p(\tau^{R*},e^{R*}) \partial_e p(\tau^{R*},e^{R*}) - \partial_\tau p(\tau^{R*},e^{R*}), \\
\label{eq:whitham_condition4}
a^R b^R &> p(\tau^R,e^R) \partial_e p(\tau^R,e^R) - \partial_\tau p(\tau^R,e^R). 
\end{align}
Moreover, we assume that the pressure law satisfies assumption \eqref{ass:convex_f}.

 Then for all $i \in \mathbb{Z}$, the updated state $\omega_i^{n+1}$, computed with the relaxation scheme \eqref{eq:relaxatation_scheme} under the CFL condition \eqref{eq:cfl_condition}, satisfies the discrete entropy inequality
 \begin{equation}
 \label{eq:discrete_entropy_inequality}
\rho_i^{n+1} \mathcal{F}(s_i^{n+1}) -
  \rho_i^{n} \mathcal{F}(s_i^{n})- \frac{ \Delta t}{\Delta x} \left( \{\rho \mathcal{F}(s) u \}_{i+1/2}^n - \{\rho \mathcal{F}(s) u \}_{i-1/2}^n \right) \leq 0,
 \end{equation}
 where we define the numerical entropy flux by
\begin{align}
\label{eq:entropy_definition_entropy_flux}
&\{\rho \mathcal F u \}_{i-1/2}^n = \{ \rho \mathcal{F}(s) u \} \left(W^{eq}(\omega_{i-1}^n),W^{eq}(\omega_i^n)\right),\\
\label{eq:entropy_entropy_flux_LR}
\{\rho \mathcal F u \}^{L,R} = &\{ \rho \mathcal{F}(s) u \} \left(W^{eq}(\omega^L),W^{eq}(\omega^R)\right) = 
\begin{cases}
      \rho^L \mathcal{F}(s(\tau^L,e^L))u^L, & \text{ if } \sigma^- > 0, \\
      \rho^{L*} \mathcal{F}(\hat{s}(W^{L*}))v^*, & \text{ if } \sigma^- < 0 \leq \sigma^v, \\
      \rho^{R*} \mathcal{F}(\hat{s}(W^{R*}))v^*, & \text{ if } \sigma^v < 0 < \sigma^{+}, \\
      \rho^R \mathcal{F}(s(\tau^R,e^R))u^R, & \text{ if } \sigma^+ < 0.
\end{cases}
\end{align}
\end{theorem}

\setcounter{remark}{1}
\begin{remark}
At the beginning of this theorem, we assume the intermediate states of density and internal energy to be positive. In Section \ref{sec:positivity} we show that the approximate Riemann solver (\ref{eq:approximate_riemann_solver}) satisfies this property for suitably chosen relaxation speeds.
\end{remark} 

\begin{proof}{\bf (Proof of Theorem \ref{theo:entropy}.)}
The proof of this theorem closely follows the steps of a similar proof in \cite[p. 113]{2}. Therefore, we only give a sketch of the proof here and do not prove every intermediate step. For more details see \cite{2}. \\
First of all, it is easy to check that 
\begin{equation}
\label{eq:entropy_riemann_invariants}
I(W) = \pi + ab \tau \quad \text{and} \quad J(W) = e-\frac{M^2(v-u)^2}{2(\frac{a}{b}-1)}-\frac{\pi^2}{2ab}
\end{equation}
are strong Riemann invariants of (\ref{sys:relaxation_model})$_{\varepsilon=\infty}$. Therefore, weak solutions of (\ref{sys:relaxation_model})$_{\varepsilon=\infty}$ satisfy 
\begin{equation}
\partial_t \rho \Psi(I,J)+\partial_x\rho\Psi(I,J)v = 0
\end{equation}
for all smooth functions $\Psi:\mathbb{R}^2 \rightarrow \mathbb{R}$.
As a consequence, for a function $W \mapsto \hat{s}(W)$, which only depends on $I$ and $J$, weak solutions of (\ref{sys:relaxation_model})$_{\varepsilon=\infty}$ satisfy the additional conservation law
\begin{equation}
\label{eq:entropy_additional_conservation_law}
\partial_t \rho \mathcal{F}(\hat{s}) + \partial_x \rho \mathcal{F}(\hat{s})v=0.
\end{equation}
We define the function $\hat{s}$ by
\begin{equation}
\label{eq:entropy_definition_hats}
\hat{s}(W) = s(\hat{\tau}(I(W),J(W)),\hat{e}(I(W),J(W))),
\end{equation}
where  $\hat{\tau}(I,J)$ is the the largest root within $\mathbb{R}^+$  of the function $f_{I,J}:\mathbb{R}^+ \rightarrow \mathbb{R}$ defined by
\begin{equation}
f_{I,J}(\tau) = \tau p\left( \tau,e(\tau,I-ab\tau)\right) + ab \tau^2 -I \tau
\end{equation}
and $\hat{e}$ is defined by
\begin{equation}
\hat{e}(I,J) = e(\hat{\tau}(I,J),I-ab \hat{\tau}(I,J)).
\end{equation}
For the further steps, the following assumption is made about the pressure law.

\setcounter{assumption}{2}
\begin{assumption}
\label{ass:convex_f}
We assume that the pressure law is such that the function $\tau \mapsto f_{I,J}$ is strictly convex for all fixed pair $(I,J)$.
\end{assumption}

\noindent
This condition is fulfilled by most common pressure laws, including the ideal gas law  \cite{2}. \\
Under this assumption, it can be proven (see \cite{2}) that for all $W$, for which the pair $(I(W),J(W))$ is in 
\begin{align}
\mathcal{A} =  \{ &(I,J)\in\mathbb{R}^2, \exists \tau >0, \exists e>0, \exists v, \exists u \ \text{such that:} \nonumber  \\ 
\label{eq:entropy_set_A_c1}
&I=p(\tau,e)+ab \tau, \\
\label{eq:entropy_set_A_c2}
&J = e-\frac{p(\tau,e)^2}{2ab}, \\
\label{eq:entropy_set_A_c3}
&ab > p(\tau,e) \partial_e p(\tau,e) - \partial_\tau p(\tau,e) \},
\end{align}
the function $\hat{s}$ is larger than the specific entropy of the original system, i.e.
\begin{equation}
\label{eq:entropy_s_inequality}
\hat{s}(W) \geq s(\tau,e) 
\end{equation} 
and that equality is reached in the relaxation equilibrium, i.e.
\begin{equation}
\label{eq:entropy_s_equilibrium}
\hat{s}(W_{|\pi=p(\tau,e),v=u}) = s(\tau,e).
\end{equation} 
Let us now go back to the additional conservation law \eqref{eq:entropy_additional_conservation_law} and integrate it over $[0,\Delta x/2)\times [0,\Delta t)$
\begin{equation}
\begin{split}
\int_0^{\Delta x/2} (\rho \mathcal{F}(\hat{s}) \left( W_\mathcal{R} \left(\frac{x}{\Delta t};W^{eq}(\omega_L),W^{eq}(\omega_R)\right) \right)
=
\int_0^{\Delta x/2} (\rho \mathcal{F}(\hat{s}))(W(x,0))dx \\
-\Delta t (\rho \mathcal{F}(\hat{s})v) \left( W_\mathcal{R} \left(\frac{\Delta x}{2 \Delta t}; W^{eq}(\omega_L),W^{eq}(\omega_R) \right) \right) \\
+ \Delta t (\rho \mathcal{F}(\hat{s})v) (W_\mathcal{R} (0;W^{eq}(\omega_L),W^{eq}(\omega_R))).
\end{split}
\end{equation}
Under consideration of the CFL condition \eqref{eq:cfl_condition} and equality \eqref{eq:entropy_s_equilibrium}, this can be rewritten as
\begin{equation}
\label{eq:entropy_intermediate_equality}
\begin{split}
&\frac{1}{\Delta x} \int_0^{\Delta x/2} (\rho \mathcal{F}(\hat{s}) \left( W_\mathcal{R} \left(\frac{x}{\Delta t};W^{eq}(\omega_L),W^{eq}(\omega_R)\right) \right) dx \\
&=
\frac{\rho^R \mathcal{F}(s^R)}{2} - \frac{\Delta t}{\Delta x} \left( \rho^R \mathcal{F}(s^R)u^R-\{\rho \mathcal F u \}^{L,R} \right). 
\end{split}
\end{equation}
The replacement of $v$ by $u$ in the entropy fluxes is due to the fact that the input values of the approximate Riemann solver are at equilibrium and therefore left and right states of $u$ and $v$ are equal in each case. Just in the intermediate states both velocities differ, which is the reason why we write $v^*$ in \eqref{eq:entropy_entropy_flux_LR}. Due to the inequality \eqref{eq:entropy_s_inequality}, it follows
\begin{equation}
\label{eq:entropy_hatsWR_sWR}
\hat{s} \left( W_\mathcal{R} \left( \frac{x}{\Delta t}; W^{eq}(\omega^L),W^{eq}(\omega^R) \right) \right) 
\geq 
s \left( (\tau^{eq},e^{eq}) \left(\frac{x}{\Delta t}; \omega^L,\omega^R \right) \right).
\end{equation}
The quantities $\tau^{eq}, e^{eq}$ on the right hand side originate from the approximate Riemann solver $W_\mathcal{R} (x / \Delta t; W^{eq}(\omega^L),W^{eq}(\omega^R))$. Since we assume $\mathcal{F}$ to be increasing, it in turn follows that
\begin{equation}
\mathcal{F}(\hat{s}) \left( W_\mathcal{R} \left( \frac{x}{\Delta t}; W^{eq}(\omega^L),W^{eq}(\omega^R) \right) \right) 
\geq
\mathcal{F}(s) \left( W_\mathcal{R}^{(\rho,\rho u, E)} \left(\frac{x}{\Delta t}; W^{eq}(\omega^L),W^{eq}(\omega^R) \right) \right).
\end{equation}
By replacing the content of the integral in \eqref{eq:entropy_intermediate_equality}, we obtain the inequality 
\begin{equation}
\begin{split}
\frac{1}{\Delta x} \int_{0}^{\Delta x/2} (\rho \mathcal{F}(s)) \left( W_\mathcal{R}^{(\rho,\rho u,E)} \left( \frac{x}{\Delta t};W^{eq}(\omega^L),W^{eq}(\omega^R) \right) \right) dx \\
\leq
\frac{\rho^R \mathcal{F}(s^R)}{2} - \frac{\Delta t}{\Delta x} \left( \rho^R \mathcal{F}(s^R) u^R -\{\rho \mathcal{F}(s) u \}^{L,R} \right).
\end{split}
\end{equation}
Inserting $\omega^L=\omega_{i-1}^n$ and $\omega^R=\omega_i^n$ leads to
\begin{equation}
\label{eq:entropy_upper_inequality}
\begin{split}
\frac{1}{\Delta x} \int_{x_{i-1/2}}^{x_i} (\rho \mathcal{F}(s)) \left( \frac{x-x_{i-1/2}}{\Delta t};\omega_{i-1}^n,\omega_i^n \right)  dx \\ 
\leq
\frac{\rho_i^n \mathcal{F}(s_i^n)}{2} - \frac{\Delta t}{\Delta x} \left( \rho_i^n \mathcal{F}(s_i^n) u_i^n -\{\rho \mathcal{F}(s) u \}_{i-1/2}^n \right).
\end{split}
\end{equation}
For the other half of the cell, on the other hand, integrating over $(-\Delta x/2,0 ] \times [0,\Delta t)$ and applying similiar steps as before results in
\begin{equation}
\begin{split}
\frac{1}{\Delta x} \int_{-\Delta x/2}^0 (\rho \mathcal{F}(s)) \left( W_\mathcal{R}^{(\rho,\rho u,E)} \left( \frac{x}{\Delta t};W^{eq}(\omega^L),W^{eq}(\omega^R) \right) \right) dx \\
\leq
\frac{\rho^L \mathcal{F}(s^L)}{2} - \frac{\Delta t}{\Delta x} \left( \{\rho \mathcal{F}(s) u \}^{L,R} - \rho^L \mathcal{F}(s^L) u^L \right),
\end{split}
\end{equation}
and inserting $\omega^L=\omega_{i}^n$ and $\omega^R=\omega_{i+1}^n$ leads to
\begin{equation}
\label{eq:entropy_lower_inequality}
\begin{split}
\frac{1}{\Delta x} \int_{x_{i}}^{x_{i+1/2}} (\rho \mathcal{F}(s)) \left( \frac{x-x_{i+1/2}}{\Delta t};\omega_{i}^n,\omega_{i+1}^n \right)  dx \\ 
\leq
\frac{\rho_i^n \mathcal{F}(s_i^n)}{2} - \frac{\Delta t}{\Delta x} \left( \{\rho \mathcal{F}(s) u \}_{i+1/2}^n - \rho_i^n \mathcal{F}(s_i^n) u_i^n \right).
\end{split}
\end{equation}
Summing up the inequalities \eqref{eq:entropy_upper_inequality} and \eqref{eq:entropy_lower_inequality} results in the inequality
\begin{equation}
\label{eq:entropy_full_inequality}
\frac{1}{\Delta x} \int_{x_{i-1/2}}^{x_{i+1/2}} (\rho \mathcal{F}(s))(\omega^n(x,t^{n+1}))dx \leq \rho_i^n \mathcal{F}(s_i^n) - \frac{\Delta t}{\Delta x} \left( \{ \rho \mathcal{F}(s)u \}_{i+1/2}^n - \{ \rho \mathcal{F}(s)u \}_{i-1/2}^n \right).
\end{equation}
Since we assume $\rho \mathcal{F}(s)$ to be strictly convex, by applying Jensen's inequality we get
\begin{equation}
\rho \mathcal{F}(s) \left(  \frac{1}{\Delta x} \int_{x_{i-1/2}}^{x_{i+1/2}} \omega^n(x,t^{n+1})dx \right) 
\leq \frac{1}{\Delta x} \int_{x_{i-1/2}}^{x_{i+1/2}} (\rho \mathcal{F}(s))(\omega^n(x,t^{n+1}))dx. 
\end{equation}
Finally, we obtain the desired discrete entropy inequality 
\begin{equation}
\rho_i^{n+1} \mathcal{F}(s_i^{n+1}) \leq \rho_i^n \mathcal{F}(s_i^n) - \frac{\Delta t}{\Delta x} \left( \{ \rho \mathcal{F}(s)u \}_{i+1/2}^n - \{ \rho \mathcal{F}(s)u \}_{i-1/2}^n \right).
\end{equation}
\hfill
\end{proof}


\subsection{Prevention of checkerboard modes}
\label{sec:checkerboards}

For asmptotic preserving methods, stationary and non-constant solutions may occur in the low-Mach regime, jumping between two different values. This behaviour can arise from the fact that the divergence or gradient of a variable is supposed to be zero in the limit equations, while the discretisation of this term allows a jumping solution. Such solutions are sometimes called checkerboards modes. Of course, it is desirable to prevent the occurrence of this unphysical phenomenon.

\setcounter{theorem}{3}
\begin{theorem}
\label{theo:checkerboard}
For the relaxation scheme, velocity and pressure are constant in space for steady periodic solutions.
\end{theorem}

\begin{proof}
The proof builds on the entropy inequality of the previous section and follows the strategy of a similar proof in \cite{1}. 
First of all, using the notations used in the entropy proof, we can write
\begin{equation}
\begin{split}
\rho_i^{n+1} \mathcal{F}(s_i^{n+1})
 \leq &
  \rho_i^{n} \mathcal{F}(s_i^{n})- \frac{ \Delta t}{\Delta x} \left( \{\rho \mathcal{F}(s) u \}_{i+1/2}^n - \{\rho \mathcal{F}(s) u \}_{i-1/2}^n \right) \\
  &=  
  \frac{1}{\Delta x} \int_{x_{i-1/2}}^{x_{i+1/2}} (\rho \mathcal{F}(\hat{s}) \left( W_\mathcal{R}(t^{n+1},x) \right) dx.
  \end{split}
\end{equation}
Additionally, by applying Jensen's inequality to the left hand side we get the following inequalities
\begin{equation}
\label{eq:checkerboard_relaxation_entropy_inequalities}
\rho_i^{n+1} \mathcal{F}(s_i^{n+1})
 \leq
  \frac{1}{\Delta x} \int_{x_{i-1/2}}^{x_{i+1/2}} \rho_i^{n+1} \mathcal{F}(s_i^{n+1}) dx 
  \leq
  \frac{1}{\Delta x} \int_{x_{i-1/2}}^{x_{i+1/2}} (\rho \mathcal{F}(\hat{s}) \left( W_\mathcal{R}(t^{n+1},x) \right) dx.
\end{equation}
We now define the left-hand side of the entropy inequality (\ref{eq:discrete_entropy_inequality}) by
\begin{equation}
\label{eq:checkerboard_D}
D_i^n:= \rho_i^{n+1} \mathcal{F}(s_i^{n+1}) -  \rho_i^{n} \mathcal{F}(s_i^{n})- \frac{ \Delta t}{\Delta x} \left( \{\rho \mathcal{F}(s) u \}_{i+1/2}^n - \{\rho \mathcal{F}(s) u \}_{i-1/2}^n \right). 
\end{equation}
For steady and space periodic solutions we then have 
\begin{equation}
\sum_i D_i^n = 0.
\end{equation}
In combination with the entropy inequality (\ref{eq:discrete_entropy_inequality}) we get
\begin{equation}
D_i^n=0 \quad \forall i.
\end{equation}
From this follows directly that all the inequalities in (\ref{eq:checkerboard_relaxation_entropy_inequalities}) are replaced by equalities and therefore the entropy is equal to the relaxation entropy
\begin{equation}
\rho_i^{n+1} \mathcal{F}(s_i^{n+1})
 =
  (\rho \mathcal{F}(\hat{s}) \left( W_\mathcal{R}(t^{n+1},x) \right).
\end{equation}
In the proof of the entropy inequality it is shown that this is just the case in the relaxation equilibrium, so only if
\begin{equation}
\pi = p(\rho,e), \ u=v, \ \tau=\frac{1}{\rho}, \ \hat{s}= s.
\end{equation}
As a consequence, the following relations apply to a single Riemann problem
\begin{equation}
\label{eq:checkerboard_states_at_single_rieman_problem}
\begin{split}
&\tau^{L*} = \frac{1}{\rho^{L*}}, \ \ \tau^{R*} = \frac{1}{\rho^{R*}}, \ \ v^{*}=u^{L*}=u^{R*}, \\
\pi^{L*}=p(\rho^{L*},&e^{L*}) = p(\rho^{L*},s^{L*}), \ \ \pi^{R*}=p(\rho^{R*},e^{R*}) = p(\rho^{R*},s^{R*}).
\end{split}
\end{equation}
Since $\tau$ is a Riemann invariant for $\sigma^-$ and $\sigma^+$, it holds
\begin{equation}
\tau^{L*} = \tau^L, \ \tau^{R*}=\tau^R.
\end{equation}
We can use this fact to gain more information about the intermediate densities
\begin{equation}
\label{eq:checkerboard_rho_states}
\begin{split}
\frac{1}{\rho^{L*}}=\tau^{L*} = \frac{1}{\rho^L} \quad \Rightarrow \quad \rho^{L*}=\rho^L ,\\
\frac{1}{\rho^{R*}}=\tau^{R*} = \frac{1}{\rho^R} \quad \Rightarrow \quad \rho^{R*}=\rho^R.
\end{split}
\end{equation}
From the explicit definition of the intermediate states in (\ref{eq:IS_rhoLs}) and (\ref{eq:IS_rhoRs}) we can deduce that
\begin{align}
\label{eq:checkerboard_rho_difference_left}
&\frac{1}{\rho^{L*}} - \frac{1}{\rho^L} = \frac{M b^R \left( v^R-v^L \right) + \pi^L-\pi^R - \bar{\rho} \left(W^L,W^R \right) \left( Z^R-Z^L \right)}{a^L \left( b^L+b^R \right)} = 0, \\
\label{eq:checkerboard_rho_difference_right}
&\frac{1}{\rho^{R*}} - \frac{1}{\rho^R} = \frac{M b^L \left( v^R-v^L \right) + \pi^R-\pi^L + \bar{\rho} \left(W^L,W^R \right) \left( Z^R-Z^L \right)}{a^R \left( b^L+b^R \right)} = 0.
\end{align}
With a look at the intermediate states $u^{L*}$ and $u^{R*}$, we see that we can use (\ref{eq:checkerboard_rho_difference_left}) and (\ref{eq:checkerboard_rho_difference_right}) to get
\begin{align}
&u^{L*} = u^L + \frac{b^L}{M} \frac{ M b^R \left(v^R-v^L \right) + \pi^L-\pi^R - \bar{\rho} \left(W^L,W^R \right) \left( Z^R-Z^L \right) }{a^L \left( b^L+b^R \right) } = u^L , \\
&u^{R*} = u^R + \frac{b^R}{M} \frac{ M b^L \left(v^L-v^R \right) + \pi^L-\pi^R + \bar{\rho} \left(W^L,W^R \right) \left( Z^R-Z^L \right)}{a^R \left( b^L+b^R \right) } = u^R.
\end{align} 
Since we are at equilibrium we can conclude that
\begin{equation}
\label{eq:checkerboard_velocity_equalities}
v^*=u^{L*}=u^{R*}=u^L=u^R=v^L=v^R.
\end{equation}
In the next part we will show that the left and the right state at the interface are equal for $\pi$. 
From the Riemann invariants in (\ref{eq:riemann_invariant_+-}) we take 
\begin{equation}
I_4^{\pm} = \frac{1}{\rho}+ \frac{\pi}{2ab}.
\end{equation}
This quantity is constant across the left and right waves in the Riemann fan which means
\begin{equation}
\begin{split}
\frac{1}{\rho^{L*}} + \frac{\pi^{L*}}{2a^Lb^L} &= \frac{1}{\rho^{L}} + \frac{\pi^{L}}{2a^Lb^L},  \\
\frac{1}{\rho^{R*}} + \frac{\pi^{R*}}{2a^Rb^R} &= \frac{1}{\rho^{R}} + \frac{\pi^{R}}{2a^Rb^R}.
\end{split}
\end{equation}
It has already been established in (\ref{eq:checkerboard_rho_states}) that the density has only two states and therefore we can simplify the equations to
\begin{equation}
\begin{split}
\pi^{L*} &= \pi^{L}  \\
\pi^{R*} &= \pi^{R}.
\end{split}
\end{equation}
From the explicit definition of the intermediate states and the closure equation
 (\ref{eq:closure_equation}) follows
\begin{align}
\pi^{L*} = \pi^L 
&= \frac{b^R \pi^L + b^L \pi^R + M b^L b^R \left(v^L-v^R \right) - b^L \bar{\rho} \left(W^L,W^R \right) \left( Z^R-Z^L \right) }{b^L+b^R} \nonumber \\
&\overset{(\ref{eq:checkerboard_velocity_equalities})}{=} \frac{b^R \pi^L + b^L \pi^R - b^L \bar{\rho} \left(W^L,W^R \right) \left( Z^R-Z^L \right) }{b^L+b^R} \nonumber \\
&\overset{(\ref{eq:closure_equation})}{=} \frac{b^R \pi^L + b^L \pi^R + b^L \left( \pi^R-\pi^L \right) }{b^L+b^R}.
\end{align}
Solving for $\pi^L$ gives
\begin{equation}
\pi^L=\pi^R.
\end{equation}
Thus we have shown that for both velocities and the pressure, the left and right states at the interface are equal. The solution in these quantities is therefore constant in space.
\hfill
\end{proof}

\setcounter{remark}{4}
\begin{remark}
\label{rem:checkerboard_isentropic}
For pressure laws that depend only on density, it can also be proven that the density and the internal energy are constant for steady and space periodic solutions.
\end{remark}

For steady periodic solutions of the relaxation method the velocity and pressure are constant, which contradicts the non-constant nature of checkerboard modes. The result of the above lemma can thus be interpreted that in velocity and pressure no checkerboard modes can occur.


\subsection{Positivity preserving property}
\label{sec:positivity}

For the robustness of a scheme it is essential to keep especially the density but also the internal energy 
positive. The following lemma will guarantee this property.\\

\setcounter{lemma}{5}
\begin{lemma}
\label{lem:positivity}
Given $\omega^L,\omega^R \in \Omega$. If the relaxation speeds $a^L$ and $a^R$ are large enough to ensure
\begin{align}
\label{eq:positivity_condition_1}
v^L - \frac{a^L}{M\rho^L} < &v^* < v^R+\frac{a^R}{M\rho^R}, \\
\label{eq:positivity_condition_2}
e^L + \frac{(\pi^{L*})^2-(\pi^L)^2}{2 a^L b^L} +& \frac{(v^*-u^{L*})^2-(v^L-u^L)^2}{2 (\frac{a^L}{b^L}-1)} >0, \\
\label{eq:positivity_condition_3}
e^R + \frac{(\pi^{R*})^2-(\pi^R)^2}{2 a^R b^R} +& \frac{(v^*-u^{R*})^2-(v^R-u^R)^2}{2 (\frac{a^R}{b^R}-1)} >0,
\end{align}
then the approximate Riemann solver $W_\mathcal{R}$ preserves the positivity of the density and internal energy.
\end{lemma}

\begin{proof}
First, it is trivial that the conditions (\ref{eq:positivity_condition_1}), (\ref{eq:positivity_condition_2}) and (\ref{eq:positivity_condition_3}) are satisfied for a sufficiently large $a$. To prove the positivity of the density in a next step, we start with the Riemann invariants $I_1^{\pm}$ from lemma \ref{lem:wave_speeds_and_riemann_invariants}, which give us
\begin{equation}
v^L - \frac{a^L}{M\rho^L} = v^* - \frac{a^L}{M \rho^{L*}} \quad \text{and} \quad v^R + \frac{a^R}{M\rho^R} = v^* + \frac{a^R}{M \rho^{R*}}. 
\end{equation}
Using these relations, we can rewrite (\ref{eq:positivity_condition_1}) by
\begin{equation}
-\rho^{L*} < 0 < \rho^{R*}.
\end{equation}
So, the intermediate states for the density are positive.
The positivity of the internal energy directly follows from (\ref{eq:positivity_condition_2}) and (\ref{eq:positivity_condition_3}), since the left-hand sides of these conditions represent the left and right intermediate states of the internal energy.
\hfill
\end{proof}\\

Clearly, this lemma is of limited use in practice. It states that in principle it is possible to preserve the positivity, but it does not help to find a suitable definition of the relaxation speeds that works generally. 
The following lemma gives stricter conditions for the relaxation speeds, which can also be used for their explicit definition. Under these conditions, it can be proven that the density is kept positive.

\begin{lemma}
\label{lem:positivity_density}
Consider the relaxation solver with intermediate values
and speeds defined by (\ref{eq:IS_vs})-(\ref{eq:IS_abZ}) with the initial data at equilibrium. 
Assume that the relaxation speeds $a^L$, $a^R$, $b^L$, $b^R$ satisfy
\begin{align}
\label{eq:positivity_speed_condition1}
a^L \geq b^L&, \quad a^R \geq b^R, \\
\label{eq:positivity_speed_condition2}
\frac{b^L}{\rho^L} \geq a_q^L&, \quad \frac{b^R}{\rho^R} \geq a_q^R, \\
\label{eq:positivity_speed_condition3}
\frac{\sqrt{a^L b^L}}{\rho^L} \geq c^L \left( 1+\beta X^L \right)&, \quad \frac{\sqrt{a^R b^R}}{\rho^R} \geq c^R \left( 1+\beta X^R \right).
\end{align}
for some $a_q^L$ and $a_q^R$ depending on $\omega^L$, $\omega^R$ and $X^L$, $X^R$ defined by (\ref{eq:positivity_XL}) and (\ref{eq:positivity_XR}) with a parameter $\beta \geq 1$. The quantities $c^L,c^R$ represent the sound speed.
Then the approximate Riemann solver $W_\mathcal{R}$ preserves the positivity of the density.
\end{lemma}

\begin{proof}
We start with the definition of the intermediate density (\ref{eq:IS_rhoLs})
\begin{align}
\frac{1}{\rho^{L*}} &= \frac{1}{\rho^L} + \frac{M b^R \left( v^R-v^L \right) + \pi^L-\pi^R - \bar{\rho} \left(W^L,W^R \right) \left( Z^R-Z^L \right)}{a^L \left( b^L+b^R \right)} \nonumber \\
&\geq \frac{1}{\rho^L} - \frac{M b^R (v^L-v^R)_+}{a^L(b^L+b^R)} -  \frac{\left(\pi^R-\pi^L+\bar{\rho}(W^L,W^R) (Z^R-Z^L)\right)_+}{a^L (b^L+b^R)} \nonumber \\
&\geq \frac{1}{\rho^L} - \frac{M (v^L-v^R)_+}{a^L} - \frac{\left(\pi^R-\pi^L+\bar{\rho}(W^L,W^R) (Z^R-Z^L)\right)_+}{a^L (\rho^L a_q^L+\rho^R a_q^R)}.
\end{align}
Analogously, we get
\begin{equation}
\frac{1}{\rho^{R*}} \geq \frac{1}{\rho^R} - \frac{M (v^L-v^R)_+}{a^R} - \frac{\left(\pi^L-\pi^R + \bar{\rho}(W^L,W^R) (Z^L-Z^R) \right)_+}{a^R (\rho^L a_q^L+\rho^R a_q^R)}.
\end{equation}
Let us now define the variables
\begin{align}
\label{eq:positivity_XL}
X^L &= \frac{1}{c^L} \left[ M \left(v^L-v^R \right)_+ + \frac{\left( \pi^R-\pi^L + \bar{\rho}(W^L,W^R) (Z^R-Z^L) \right)_+}{\rho^L a_q^L+\rho^R a_q^R} \right], \\
\label{eq:positivity_XR}
X^R &= \frac{1}{c^R} \left[ M \left(v^L-v^R \right)_+ + \frac{\left( \pi^L-\pi^R + \bar{\rho}(W^L,W^R) (Z^L-Z^R) \right)_+}{\rho^L a_q^L+\rho^R a_q^R} \right],
\end{align}
in order to rewrite the former inequalities in the form
\begin{equation}
\label{eq:positivity_inequality_rhos}
\begin{split}
\frac{1}{\rho^{L*}} &\geq  \frac{1}{\rho^L} \left( 1-\frac{\rho^L c^L}{a^L} X^L \right), \\
\frac{1}{\rho^{R*}} &\geq \frac{1}{\rho^R} \left( 1-\frac{\rho^R c^R}{a^R} X^R \right).
\end{split}
\end{equation}
From combining the conditions (\ref{eq:positivity_speed_condition1}) and (\ref{eq:positivity_speed_condition3}), it follows that
\begin{equation}
\begin{split}
\frac{a^L}{\rho^L} &\geq c^L (1+\beta X^L) \Rightarrow \frac{\rho^L c^L}{a^L}	\leq \frac{1}{1+\beta X^L}, \\
\frac{a^R}{\rho^R} &\geq c^R (1+\beta X^R) \Rightarrow \frac{\rho^R c^R}{a^R}	\leq \frac{1}{1+\beta X^R}.
\end{split}
\end{equation}
With these inequalities we rewrite (\ref{eq:positivity_inequality_rhos})
\begin{equation}
\begin{split}
\frac{1}{\rho^{L*}} &\geq  \frac{1}{\rho^L} \left( 1-\frac{X^L}{1+\beta X^L} \right), \\
\frac{1}{\rho^{R*}} &\geq \frac{1}{\rho^R} \left( 1-\frac{X^R}{1+\beta X^R} \right).
\end{split}
\end{equation}
Because of the definitions in (\ref{eq:positivity_XL}) and (\ref{eq:positivity_XR}) we know that $X^L,X^R\geq 0$ and therefore we can conclude that
\begin{equation}
\rho^{L*} > 0, \quad \rho^{R*} > 0.
\end{equation}
\hfill
\end{proof}\\

A similar proof for the positivity of the internal energy would be complicated due to the more complex structure of its intermediate states. Therefore, we choose a different way of proof, based on the proof of the entropy inequality in Sect. \ref{sec:entropy}.

\begin{lemma}
Under the conditions of the entropy theorem \ref{theo:entropy} and for given $\omega^L,\omega^R \in \Omega$, the relaxation scheme using the approximate Riemann solver $W_\mathcal{R}$ defined in \eqref{eq:approximate_riemann_solver} preserves the positivity of the internal energy.
\end{lemma}

\begin{proof}
In the proof of the entropy inequality, it is stated in \eqref{eq:entropy_s_inequality} that the specific relaxation entropy is larger than the specific entropy of the original system. Therefore, we can conclude
\begin{equation}
 \hat{s}_i^{n+1} \geq s(\rho_i^{n+1},e_i^{n+1})=s_i^{n+1}.
\end{equation}
In a next step we show that $\hat{s}_i^{n+1}$ is positive. For this we consider one Riemann problem and investigate the input of the function $\hat{s}(W)=s(\hat{\tau}(I(W),J(W)),\hat{e}(I(W),J(W)))$. Here $\hat{\tau}$ is already positive by definition. For $\hat{e}(I(W),J(W))$ depending on $W=W^L$ or $W=W^R$ its positivity is trivial. For $W=W^{L*}$ we can rewrite $\hat{e}$ using $J$ and additionally make use of the fact that $I$ and $J$ are strong Riemann invariants for $\sigma^L$. It follows
\begin{align*}
\hat{e}^{L*} &=  J^{L} + \frac{M^2(v^*-u^{L*})^2}{2(\frac{a}{b}-1)} + \frac{(I^{L}-a^Lb^L\hat{\tau}^{L*})^2}{2a^Lb^L} \\
&= e^L - \frac{M^2(v^L-u^L)^2}{2(\frac{a^L}{b^L}-1)} - \frac{(\pi^L)^2}{2a^Lb^L} + \frac{M^2(v^*-u^{L*})^2}{2(\frac{a^L}{b^L}-1)} + \frac{(\pi^L+a^Lb^L\tau^L-a^Lb^L\hat{\tau}^{L*})^2}{2a^Lb^L} \\
&\geq e^L + \frac{M^2(v^*-u^{L*})^2}{2(\frac{a^L}{b^L}-1)} + \frac{(a^L b^L \tau^L - a^L b^L \hat{\tau}^{L*})^2}{2 a^L b^L}.
\end{align*}
Then for given $e^L>0$ it follows that $\hat{e}^{L*}>0$. The same arguments lead to $\hat{e}^{R*}>0$. The positivity of $\hat{s}$ then follows from the definition of the function $s(\tau,e): \mathbb{R}^+ \times \mathbb{R}^+ \rightarrow \mathbb{R^+}$.\\
Furthermore, from relation (\ref{eq:temperature}) it follows that
\begin{equation}
\partial_s e(\tau,s) = - T(\tau,e) < 0.
\end{equation}
For the positive specific relaxation entropy $\hat{s}_i^{n+1}$ and the definition of $e(\rho,s):\mathbb{R}^+\times \mathbb{R}^+ \rightarrow \mathbb{R}^+$ we then obtain a positive updated internal energy
\begin{equation}
e_i^{n+1} = e(\rho_i^{n+1},s_i^{n+1}) \geq e \left(\rho_i^{n+1},\hat{s}_i^{n+1} \right) > 0.
\end{equation}
\hfill
\end{proof}


\subsection{Asymptotic preserving property}
\label{sec:AP}

In the low Mach limit, the solutions of the Euler equations (\ref{sys:Euler}) tend to the solutions of the incompressible Euler equations. This behaviour can be illustrated by inserting expansions in terms of $M$ given by
\begin{equation}
\label{eq:ap_expansions_in_M}
\begin{split}
\rho &= \rho_0 + M \rho_1 + M^2 \rho_2 + \mathcal{O}(M^3), \quad \boldsymbol{u} = \boldsymbol{u}_0 + M \boldsymbol{u}_1 + M^2 \boldsymbol{u}_2 + \mathcal{O}(M^3), \\
e &= e_0+ M e_1 + M^2 e_2 + \mathcal{O}(M^3), \quad \ p = p_0 + M p_1 + M^2 p_2 + \mathcal{O}(M^3),
\end{split}
\end{equation}
into the Euler equations (\ref{sys:Euler}). Now one can collect terms of order $\mathcal{O}(M^{-2})$
\begin{equation}
\label{eq:ap_order_M-2_terms}
\nabla p_0 = -\rho_0 \nabla \Phi,
\end{equation}
of order $\mathcal{O}(M^{-1})$
\begin{equation}
\label{eq:ap_order_M-1_terms}
\nabla p_1 = -\rho_1 \nabla \Phi,
\end{equation}
and finally of order $\mathcal{O}(1)$
\begin{equation}
\label{eq:ap_order_M0_terms}
\begin{split}
\nabla \cdot \left( \rho_0 \boldsymbol{u}_0 \right) &= 0, \\
\partial_t \boldsymbol{u}_0 + \boldsymbol{u}_0 \cdot \nabla \boldsymbol{u}_0 + \frac{\nabla p_2}{\rho_0} &= -\frac{\rho_2 \nabla \Phi}{\rho_0},  \\
\partial_t e_0 + \boldsymbol{u}_0 \cdot \nabla e_0 + \frac{1}{\rho_0} \nabla \cdot (p_0 \boldsymbol{u}_0 ) &= - \boldsymbol{u}_0 \cdot \nabla \Phi.
\end{split}
\end{equation}
These equations describe incompressible flows. The conditions (\ref{eq:ap_order_M-2_terms}) and (\ref{eq:ap_order_M-1_terms}) show that the couples $\rho_0,p_0$ and $\rho_1,p_1$ fulfil the hydrostatic equilibrium and are therefore time independent. This property is used in the derivation of the limit equations in (\ref{eq:ap_order_M0_terms}).

\setcounter{remark}{8}
\begin{remark}
The third equation in (\ref{eq:ap_order_M0_terms}) vanishes because of (\ref{eq:ap_order_M-2_terms}). In consequence, the limit equations contain an unknown  $\rho_2$, 
to which no conditions seem to be attached that determine its behavior. Under the assumption that the term $\frac{\nabla p_2}{\rho_0} + \frac{\rho_2 \nabla \Phi}{\rho_0}$ can be written as a gradient, this term can be replaced by $\nabla \Pi$, so that the second equation in (\ref{eq:ap_order_M0_terms}) changes to
\begin{equation}
\partial_t \boldsymbol{u} + \boldsymbol{u} \cdot \nabla \boldsymbol{u} + \nabla \Pi = 0.
\end{equation}
With this equation all variables can be determined. For similar arguments, see \cite{35}.
\end{remark}

In the next step, we want to analyse to what extent the solutions of the compressible Euler equations correspond to those of the incompressible equations. Under the assumptions that in the density no constant fluctuations occur, i.e.
\begin{equation}
\rho = \rho_0 + \mathcal{O}(M^2),
\end{equation}
and that the hydrostatic equilibrium is fulfilled up to errors of order $\mathcal{O}(M^2)$
\begin{equation}
\nabla p + \rho \nabla \Phi = \mathcal{O}(M^2),
\end{equation}
the Euler equations (\ref{sys:Euler}) become
\begin{equation}
\label{sys:ap_limit_Euler_equations}
\begin{split}
\nabla \cdot (\rho \boldsymbol{u}) &= \mathcal{O} (M^2), \\
\partial_t \boldsymbol{u} + \boldsymbol{u} \cdot \nabla \boldsymbol{u} + \frac{\nabla p_2}{ \rho_0} &= -\frac{\rho_2 \nabla \Phi}{ \rho_0} + \mathcal{O}(M^2), \\
\partial_t e + \boldsymbol{u} \cdot \nabla e + \frac{1}{\rho_0} \nabla \cdot (p \boldsymbol{u}) &= - \boldsymbol{u} \cdot \nabla \Phi + \mathcal{O}(M^2).
\end{split}
\end{equation}
The solutions of (\ref{sys:Euler}) thus agree with those of the incompressible model up to an error of order $M^2$. \\

Following these theoretical results, the numerical scheme should be consistent with the limit behaviour as M tends to zero, in the sense that the discretization for the compressible Euler equations should tend to the incompressible Euler equations when the Mach number tends to zero. 
The key to achieve this behaviour for the presented relaxation scheme is the definition of the relaxation speeds $a$ and $b$. 
In the former sections several conditions are imposed on these speeds that have to be satisfied so that the scheme is stable and has the properties presented in the former sections. A suitable choice that indeed fulfils the so far stated requirements is the classical one, in which $a$ and $b$ are set to be equal
\begin{equation}
\begin{split}
\label{eq:ap_relaxation_speeds_without_fix}
a_q^\alpha &= c^\alpha, \\
a^\alpha &= b^\alpha = \rho^\alpha c^\alpha (1+\beta X^\alpha).
\end{split}
\end{equation}
This definition closely follows the condition (\ref{eq:positivity_speed_condition3}) in Lemma \ref{lem:positivity_density}. Unfortunately, this definition does not lead to an appropriate discretization, but to excessive diffusion in the low Mach limit.

In order to change this behaviour the speeds have to be redefined. In this context it is important to ensure that not only the diffusion is reduced, but also that the sub-characteristic condition remains fulfilled.
A suitable choice proposed by the authors of \cite{1} is given by
\begin{equation}
\begin{split}
\label{eq:ap_relaxation_speeds_with_fix}
a_q^\alpha &= \min (1,M) c^\alpha, \\
a^\alpha &= \frac{\rho^\alpha}{\min (1,M)} c^\alpha (1+\beta X^\alpha), \\
b^\alpha &=  \min (1,M) \rho^\alpha c^\alpha (1+\beta X^\alpha).
\end{split}
\end{equation}
By this definition the speeds are rescaled in the case of small Mach numbers, i.e. for $M< 1$.

\begin{remark}
\label{rem:relaxation_speeds}
In the case of Mach numbers $M \geq 1$, the relaxation speeds are equal ($a=b$) and we obtain a classical relaxation system with only one relaxation speed.
\end{remark}

\begin{remark}
\label{rem:CFL_condition}
The new scaling of the relaxation speed $a$ has the effect that the maximum wave speed increases by an order of magnitude $M$. As a consequence, the CFL condition (\ref{eq:cfl_condition}) becomes stricter and the time step must be chosen smaller accordingly, i.e.
\begin{equation}
\Delta t \sim \frac{M^2 \Delta x}{c}.
\end{equation}
As shown in \cite{1}, by replacing $M$ by $\hat{M}=\max \{M^2,k\Delta x \}$ in the relaxation scheme the CFL condition can be reduced to the parabolic-type condition
\begin{equation}
\Delta t \sim \frac{\max \{M^2,k\Delta x \} \Delta x}{c}.
\end{equation}
\end{remark}

\setcounter{theorem}{11}
\begin{theorem}
\label{theo:asymptotic_preserving}
The two-speed relaxation scheme with the relaxation speeds (\ref{eq:ap_relaxation_speeds_with_fix}) is asymptotic preserving in the sense that:
\begin{enumerate}[a)]
\item it is first order uniformly with respect to the Mach number M and
\item for $M<\sqrt{k \Delta x}$ and $k$ constant it is consistent at first order with the incompressible limit model (\ref{eq:ap_order_M0_terms}).
\end{enumerate}
\end{theorem}

\begin{proof}
In order to prove the first statement of the theorem we evaluate the consistency error by expanding the numerical flux (\ref{eq:numerical_flux}) in terms of $M$ and then subtract the central flux $(F( \omega^L)+F(\omega^R))/2$.

In the low Mach limit $M \rightarrow 0$, the wave speeds $\sigma^-$ and $\sigma^+$ in (\ref{eq:wave_speeds}) tend towards infinity. Therefore it is sufficient just to consider the intermediate fluxes $F^{L*}$ and $F^{R*}$ for the numerical flux. 
In a first step of the analysis we rewrite the relaxation speeds as expansions in terms of $M$, so we get
\begin{equation}
X^{\alpha} = \mathcal{O}(M), \quad \  b^{\alpha} = M \bar{b}^\alpha + \mathcal{O}(M^2), \quad \ a^{\alpha} = \frac{\bar{b}^{\alpha}}{M} (1+\mathcal{O}(M))
\end{equation}
with
\begin{equation}
\bar{b}^{\alpha} = \rho^{\alpha} c^{\alpha}.
\end{equation}
Since
\begin{equation}
\bar{b}^{R}-\bar{b}^{L} = \mathcal{O}(M^2),
\end{equation}
we can write $\bar{b}$ instead of $\bar{b}^{L}$ and $\bar{b}^{R}$ up to errors of $\mathcal{O}(M^2)$. Expanding the intermediate states (\ref{eq:IS_vs})-(\ref{eq:IS_piRs}) in terms of $M$ yields
\begingroup
\allowdisplaybreaks
\begin{align}
v^* =& \frac{u^L+u^R}{2} + \frac{\pi^L-\pi^R}{2 M^2 \bar{b}}-\frac{\rho (Z^R-Z^L)}{2 M^2 \bar{b}}\\
 &+ \mathcal{O}(M(u^L-u^R)) + \mathcal{O}( \frac{\pi^L-\pi^R+\bar{\rho} (Z^R-Z^L)}{M} ), \\
\pi^{L*} =& \frac{\pi^L+\pi^R}{2} + M^2 \bar{b} \frac{u^L-u^R}{2} + \frac{\bar{\rho}(Z^R-Z^L)}{2\bar{b}} \\
 &+ \mathcal{O}(M^3(u^L-u^R)) + \mathcal{O}( M(\pi^L-\pi^R+\bar{\rho} (Z^R-Z^L)) ), \\
\pi^{R*} =& \frac{\pi^L+\pi^R}{2} + M^2 \bar{b} \frac{u^L-u^R}{2} - \frac{\bar{\rho}(Z^R-Z^L)}{2\bar{b}} \\
 &+ \mathcal{O}(M^3(u^L-u^R)) + \mathcal{O}( M(\pi^L-\pi^R+\bar{\rho} (Z^R-Z^L)) ), \\
\frac{1}{\rho^{L*}} =& \frac{1}{\rho^L}+ \mathcal{O}(M^2(u^L-u^R)) + \mathcal{O}( \pi^L-\pi^R+\bar{\rho} (Z^R-Z^L) ), \\
\frac{1}{\rho^{R*}} =& \frac{1}{\rho^R}+ \mathcal{O}(M^2(u^L-u^R)) + \mathcal{O}( \pi^L-\pi^R+\bar{\rho} (Z^R-Z^L) ), \\
u^{L*} =& u^L+ \mathcal{O}(M^2(u^L-u^R)) + \mathcal{O}( \pi^L-\pi^R+\bar{\rho} (Z^R-Z^L) ), \\
u^{R*} =& u^R+ \mathcal{O}(M^2(u^L-u^R)) + \mathcal{O}( \pi^L-\pi^R+\bar{\rho} (Z^R-Z^L) ).
\end{align}
\endgroup
We can derive these expansion and put the terms $\pi^L-\pi^R+\bar{\rho}(Z^R-Z^L)$ into the error estimates, since, as stated in (\ref{eq:ap_order_M-2_terms}) and (\ref{eq:ap_order_M-1_terms}), the hydrostatic equilibrium is satisfied up to terms of order $\mathcal{O}(M^2)$ in the low Mach limit, i.e.
\begin{equation}
\label{eq:ap_hydrostatic_equilibrium_M^2}
p^L-p^R + \bar{\rho} (Z^L-Z^R) = \mathcal{O}(M^2).
\end{equation}

\noindent
With the help of these expansions, we calculate the flux differences component by component.
\begin{enumerate}[i)]
\item The difference for the left intermediate flux $F^{L*}$ in the first component writes
\begin{align*}
&\rho^{L*} v^* - \frac{\rho^L u^L + \rho^R u^R}{2} \\
=& -\frac{\rho^L u^L +\rho^Ru^R}{2}+\frac{\rho^L}{2\bar{b}} \left( \frac{p^L-p^R}{M^2}+\frac{\bar{\rho}(Z^L-Z^R)}{M^2}\right) \\
&+ \rho^L \frac{u^L+u^R}{2} + \mathcal{O}(M(u^L-u^R))+ \mathcal{O}(\frac{\pi^L-\pi^R+\bar{\rho} (Z^R-Z^L)}{M}).
\end{align*}
This difference can be further simplified. In the low Mach limit, the density is constant up to errors of $\mathcal{O}(M^2)$. Therefore we can write 
\begin{equation}
\rho^R - \rho^L = \mathcal{O}(M^2)
\end{equation}
and replace $\rho^R$ in the difference by $\rho^L$. 
Additionally, we replace the differences between the left and right states by numerical derivatives, i.e.
\begin{equation}
\label{eq:ap_numerical_derivative}
\begin{split}
u^L-u^R &= -\Delta x \partial_x u + \mathcal{O}(\Delta x^2), \\
p^L-p^R &= -\Delta x \partial_x p + \mathcal{O}(\Delta x^2), \\
Z^L-Z^R &= -\Delta x \partial_x Z +\mathcal{O}(\Delta x^2). 
\end{split}
\end{equation}
Applying these simplifications results in
\begin{equation}
\rho^{L*} v^* - \frac{\rho^L u^L + \rho^R u^R}{2}
= -\frac{\Delta x}{2} \frac{\rho^L}{\bar{b}} \left( \partial_x \frac{p}{M^2} + \bar{\rho} \partial_x \frac{Z}{M^2} \right) +\mathcal{O}(\Delta x^2) + \mathcal{O}(M \Delta x).
\end{equation}
The denominator $M^2$ does not lead to excessive diffusion at this point, as again the hydrostatic equilibrium is fulfilled up to $\mathcal{O}(M^2)$. Analogous calculations for the right intermediate flux $F^{R*}$ lead to
\begin{equation}
\rho^{R*} v^* - \frac{\rho^L u^L + \rho^R u^R}{2}
= -\frac{\Delta x}{2} \frac{\rho^R}{\bar{b}} \left( \partial_x \frac{p}{M^2} + \bar{\rho} \partial_x \frac{Z}{M^2} \right) +\mathcal{O}(\Delta x^2) + \mathcal{O}(M \Delta x).
\end{equation}
\item The second component for the left flux can be expressed by
\setlength{\jot}{4pt}
\begin{align*}
&\rho^{L*} u^{L*} v^* + \frac{\pi^{L*}}{M^2} - \frac{\rho^L (u^L)^2+\frac{\pi^L}{M^2}+\rho^R (u^R)^2+\frac{\pi^R}{M^2}}{2} \\
=&\bar{b} \frac{u^L-u^R}{2} 
+ \rho^L u^L \frac{u^L+u^R}{2} 
- \rho^L u^R \frac{u^L-u^R}{2}+\rho^L u^R \frac{u^L-u^R}{2} \\
&- \frac{\rho^L(u^L)^2+\rho^R(u^R)^2}{2}
+\rho^L u^L \frac{p^L-p^R+\bar{\rho}(Z^L-Z^R)}{2\bar{b}M^2} \\
&-\frac{\bar{\rho}(Z^L-Z^R)}{2 M^2}
+ \mathcal{O}(M(u^L-u^R))+ \mathcal{O}(\frac{\pi^L-\pi^R+\bar{\rho} (Z^R-Z^L)}{M}) \\
=& \bar{b} \frac{u^L-u^R}{2} 
+\rho^L u^R \frac{u^L-u^R}{2} 
+\rho^L u^L \frac{p^L-p^R+\bar{\rho}(Z^L-Z^R)}{2\bar{b}M^2} \\
&-\frac{\bar{\rho}(Z^L-Z^R)}{2 M^2}
+ \mathcal{O}(M(u^L-u^R))+ \mathcal{O}(\frac{\pi^L-\pi^R+\bar{\rho} (Z^R-Z^L)}{M}) \\
=& - \frac{\Delta x}{2} \left( \bar{b}+\rho^L u^R \right) \partial_x u
- \frac{\Delta x}{2} \frac{\rho^L u^L}{\bar{b}} \left( \partial_x \frac{p}{M^2} + \bar{\rho} \partial_x \frac{Z}{M^2} \right) + \frac{ \Delta x}{2} \bar{\rho} \partial_x \frac{Z}{M^2} \\ &+\mathcal{O}(\Delta x^2)+ \mathcal{O}(M \Delta x)
\end{align*}
and for the right flux by
\begin{align*}
&\rho^{R*} u^{R*} v^* + \frac{\pi^{R*}}{M^2} - \frac{\rho^L (u^L)^2+\frac{\pi^L}{M^2}+\rho^R (u^R)^2+\frac{\pi^R}{M^2}}{2} \\
=& - \frac{\Delta x}{2} \left( \bar{b}+\rho^R u^L \right) \partial_x u
- \frac{\Delta x}{2} \frac{\rho^R u^R}{\bar{b}} \left( \partial_x \frac{p}{M^2} + \bar{\rho} \partial_x \frac{Z}{M^2} \right) - \frac{ \Delta x}{2} \bar{\rho} \partial_x \frac{Z}{M^2} \\ &+\mathcal{O}(\Delta x^2) + \mathcal{O}(M \Delta x).
\end{align*}
In this flux difference, the new scaling of the relaxation speeds defined in (\ref{eq:ap_relaxation_speeds_with_fix}) unfolds its importance. Clearly, the viscosity on the velocity, represented by the first term, is independent of the Mach number and therefore does not increase in the low Mach limit. With the classical scaling (\ref{eq:ap_relaxation_speeds_without_fix}), on the other hand, this term would have the size $\mathcal{O}(1/M)$ leading to excessive diffusion for low Mach numbers.
While a Mach number dependence in the first term would be problematic, it is not in the second term due to (\ref{eq:ap_hydrostatic_equilibrium_M^2}). 
The remaining third term containing the derivative of the gravitational potential, which also depends on $1/M^2$, cancels out with the gravitational source term (\ref{eq:source_terms}) in the relaxation scheme.
\item For the difference in the third component, similar steps for the left flux result in
\begin{align*}
&\left( \left( \frac{1}{2}M^2 \rho^{L*} (u^{L*})^2 + \rho^{L*} e^{L*} \right) + \pi^{L*} \right) v^*
-\frac{
 ( 
E^L
 + p^{L} 
  ) u^L 
 + ( 
E^R
  + p^{R} 
   ) u^R}{2}\\
=& \rho^L u^R \frac{e^L -e^R}{2}+u^R \frac{p^L-p^R}{2} +\frac{ \rho^L e^L+p^L }{2 \bar{b}}  \frac{p^L-p^R+\bar{\rho} (Z^L-Z^R)}{M^2} \\
&+ \mathcal{O}(M(u^L-u^R))+ \mathcal{O}(\frac{\pi^L-\pi^R+\bar{\rho} (Z^R-Z^L)}{M})\\
=& -\frac{\Delta x}{2} \rho^L u^R \partial_x e -\frac{\Delta x}{2} u^R \partial_x p - \frac{\Delta x}{2} \frac{ \rho^L e^L+p^L}{\bar{b}} \left( \partial_x \frac{p}{M^2} + \partial \frac{Z}{M^2} \right) \\
 &+\mathcal{O}(\Delta x^2) + \mathcal{O}(M \Delta x).
\end{align*}
and for the right flux in
\begin{align*}
&\left( \left( \frac{1}{2}M^2 \rho^{R*} (u^{R*})^2 + \rho^{R*} e^{R*} \right) + \pi^{R*} \right) v^*
-\frac{
 ( 
E^L
 + p^{L} 
  ) u^L 
 + ( 
E^R
  + p^{R} 
   ) u^R}{2}\\
=& \frac{\Delta x}{2} \rho^R u^L \partial_x e +\frac{\Delta x}{2} u^L \partial_x p - \frac{\Delta x}{2} \frac{ \rho^R e^R+p^R}{\bar{b}} \left( \partial_x \frac{p}{M^2} + \partial \frac{Z}{M^2} \right) \\
&+\mathcal{O}(\Delta x^2) + \mathcal{O}(M \Delta x).
\end{align*}
\setlength{\jot}{0pt}
\end{enumerate}
The expansions for all three components are first-order uniformly in $M$. It is particularly important that the viscosity on the velocity $u$ is independent of $M$. \\

The result of the first statement can now be used to prove the second statement of the theorem. We have proven that the solution $w_{M,\Delta x}$ of the relaxation scheme is consistent with the exact solution $w_M$ of the dimensionless Euler equations (\ref{sys:Euler}) up to order $\mathcal{O}(\Delta x)$ independent of the Mach number, i.e.
\begin{equation}
\label{eq:ap_consistence_relaxation_Euler}
w_{M,\Delta x} - w_M = \mathcal{O}(\Delta x).
\end{equation}
 Additionally, we can deduce from the system (\ref{sys:ap_limit_Euler_equations}) that $w_M$ is consistent with the solution $w$ of the incompressible Euler equations up to order $\mathcal{O}(M^2)$, i.e.
\begin{equation}
\label{eq:ap_consistence_Euler_incompressibleEuler}
w_{M} - w = \mathcal{O}(M^2).
\end{equation}
 Combining (\ref{eq:ap_consistence_relaxation_Euler}) and (\ref{eq:ap_consistence_Euler_incompressibleEuler}) with the condition that $M^2=\mathcal{O}(\Delta x)$ finally results in
 \begin{equation}
 \label{eq:ap_consistence_relaxation_incompressibleEuler}
 w_{M,\Delta x}-w=\mathcal{O}(\Delta x)
 \end{equation}
 and therefore meets the second statement of the theorem.
 \hfill
\end{proof}


\subsection{Well-balanced property}
\label{sec:wb}

As explained in the introduction, the well-balanced property is important for solving problems close to hydrostatic equilibrium. 
In a first step, we will show that the approximate Riemann solver satisfies this property.
Building on this, we will then prove in the second step that the entire scheme has this property.

\setcounter{lemma}{12}
\begin{lemma}
\label{lem:wb_discrete_steady_states}
Assume two given states at equilibrium $W^L$ and $W^R$ satisfy
\begin{align}
\label{eq:wb_steady_state_condition1}
&u^L = u^R = 0, \\
\label{eq:wb_steady_state_condition2}
&p^R-p^L+\bar{\rho}(W^L,W^R) (\Phi^R-\Phi^L) = 0.
\end{align}
Then the approximate Riemann solver $W_{\mathcal{R}}$ preserves the steady state, i.e.
\begin{equation}
W_{\mathcal{R}} (x/t, W^L,W^R) =
\left\{\begin{array}{ll} W^L \quad \text{if} \ x/t <0,  \\
							  W^R \quad \text{if} \ x/t>0.
										\end{array}\right.
\end{equation}
\end{lemma}

\begin{proof}
The result directly follows from the definition of the intermediate states given in (\ref{eq:IS_vs})-(\ref{eq:IS_abZ}). Consider the intermediate state $v^*$. Since we start at equilibrium, we can replace the relaxation variables by their corresponding original variables. Using the conditions (\ref{eq:wb_steady_state_condition1})-(\ref{eq:wb_steady_state_condition2}) results in
\begin{equation*}
v^* = \frac{1}{b^L+b^R}\left( M b^L u^L + M b^R u^R + p^L-p^R - \bar{\rho} \left(W^L,W^R \right) \left( \Phi^R-\Phi^L \right) \right) = 0. \\
\end{equation*}
Similar calculations for the other intermediate states complete the proof.
\hfill
\end{proof}\\

Lemma \ref{lem:wb_discrete_steady_states} is rather general, as it assumes that the conditions in (\ref{eq:wb_steady_state_condition1}) and (\ref{eq:wb_steady_state_condition2}) are satisfied.  Clearly, these conditions depend on the definition of the $\bar{\rho}$-function. For a simple definition like the arithmetic mean, which is not adjusted to the underlying hydrostatic equilibrium, the scheme maintains the equilibrium to second order \cite{2}. Since we are free to define $\bar{\rho}$ we can adjust it to the hydrostatic equilibrium and maintain it even up to machine precision. The only limiting requirement for $\bar{\rho}$ that has to be considered is the consistency property 
\begin{equation}
\label{eq:wb_rhobar_consistency}
\rho^L = \rho^R=\rho \quad  \Rightarrow \quad \bar{\rho}(W^L,W^R) = \rho. \\ 
\end{equation}
 The following lemma describes the adjusted definitions for isothermal, incompressible and polytropic equilibria. These definitions were already described in \cite{2}.
 
\begin{lemma}
\label{lem:wb_3_families_hydrostatic_equilibria}
\begin{enumerate} [i)]
\item Let $W^L$ and $W^R$ be two states satisfying the isothermal equilibrium
\begin{equation}
\label{eq:wb_isothermal_equilibrium}
\left\{\begin{array}{ll} 
u^L = u^R = 0, \\
\rho^{L,R} = \exp{ \frac{C-\Phi^{L,R}}{K}}, \\
p^{L,R} = K \exp{ \frac{C-\Phi^{L,R}}{K}},
\end{array}\right.
\end{equation}
with $K>0$ and $C\in \mathbb{R}$. If the function $\bar{\rho}$ is defined by
\begin{equation}
\label{eq:wb_rhobar_isothermal}
\bar{\rho} (W^L,W^R) = 
\left\{\begin{array}{ll} 
\frac{\rho^R-\rho^L}{\ln (\rho^R) - \ln (\rho^L)} \quad &\text{if} \ \rho^L \neq \rho^R, \\
\rho^L \quad &\text{if} \ \rho^L=\rho^R,
\end{array}\right.
\end{equation}
then the approximate Riemann solver $W_{\mathcal{R}}$ preserves the steady state.
\item Let $W^L$ and $W^R$ be two states satisfying the incompressible equilibrium
\begin{equation}
\left\{\begin{array}{ll} 
u^L = u^R = 0, \\
\rho^{L} = \rho^{R}, \\
p^{L}+\rho^L \Phi^L = p^R+\rho^R \Phi^R.
\end{array}\right.
\end{equation}
If the function $\bar{\rho}$ satisfies the consistency condition (\ref{eq:wb_rhobar_consistency}), then the approximate Riemann solver $W_{\mathcal{R}}$ preserves the steady state.
\item Let $W^L$ and $W^R$ be two states satisfying the polytropic equilibrium
\begin{equation}
\left\{\begin{array}{ll} 
u^L = u^R = 0, \\
\rho^{L,R} = \bigg ( \frac{\Gamma - 1}{\Gamma K} ( C-\Phi^{L,R} ) \bigg )^{\frac{\Gamma}{\Gamma -1}}, \\
p^{L,R} = K^{\frac{1}{1-\Gamma}} \bigg ( \frac{\Gamma - 1}{\Gamma} ( C-\Phi^{L,R} ) \bigg )^{\frac{\Gamma}{\Gamma -1}},
\end{array}\right.
\end{equation}
with $\Gamma \in (0,1) \cup (1,+\infty)$, $K>0$ and $C\in \mathbb{R}$. If the function $\bar{\rho}$ is defined by
\begin{equation}
\label{eq:wb_rhobar_polytropic}
\bar{\rho} (W^L,W^R) = 
\left\{\begin{array}{ll} 
\frac{\Gamma -1}{\Gamma} \frac{(\rho^R)^\Gamma-(\rho^L)^\Gamma}{(\rho^R)^{\Gamma-1}-(\rho^L)^{\Gamma-1}} \quad &\text{if} \ \rho^L \neq \rho^R, \\
\rho^L \quad &\text{if} \ \rho^L=\rho^R,
\end{array}\right.
\end{equation}
then the approximate Riemann solver $W_{\mathcal{R}}$ preserves the steady state.
\end{enumerate}
\end{lemma}

\begin{proof}
In order to prove this lemma it is sufficient to show that with the explicit definition of $\bar{\rho}$ the conditions (\ref{eq:wb_steady_state_condition1}) and (\ref{eq:wb_steady_state_condition2}) are satisfied. If so, we can use lemma \ref{lem:wb_discrete_steady_states} and the proof is complete.
Using the definitions of the isothermal equilibrium states, we can determine the following differences
\begin{align*}
\Phi^R-\Phi^L &= K ( \ln(\rho^R)-\ln(\rho^L)), \\
p^R-p^L &= K (\rho^R-\rho^L).
\end{align*}
By inserting these differences together with $\bar{\rho}$ defined by (\ref{eq:wb_rhobar_isothermal}) into equation (\ref{eq:wb_steady_state_condition2}), it becomes clear that this condition is satisfied. Together with the velocities, which are zero, lemma \ref{lem:wb_discrete_steady_states} can be applied and the proof of i) is complete. The proofs for incompressible and polytropic equilibria work in the same way. For more details we refer the reader to \cite{2}.
\hfill
\end{proof}

\setcounter{remark}{14}
\begin{remark}
\label{rem:wb_implementation}
To ensure the exact preservation of steady states at rest, it is important to consider the following two points in the implementation:
\begin{enumerate}
\item The comparative operators in the approximate Riemann solver (\ref{eq:numerical_flux}) must be adjusted to the definition of the sign function in the programming language used. The choice provided here is adapted to $sign(0)=1$.
\item The implementation of the classical definition of the logarithmic mean can lead to problems if left and right input are very close. Ismail and Roe provide an alternative way of implementation in \cite{26}, which avoids this problem.
\end{enumerate}
\end{remark}

In practical applications, e.g. in astrophysics, the hydrostatic states are often just available as discrete data generated by previously performed simulations. The following lemma provides an approach to maintain these hydrostatic equilibria as well.\\

\setcounter{lemma}{15}
\begin{lemma}
\label{lem:wb_arbitrary_hydrostatic_equilibrium}
Let $W^L$ and $W^R$ be two states satisfying some hydrostatic equilibrium
\begin{equation}
\label{eq:wb_arbitrary_hydrostatic_equilibrium}
\left\{\begin{array}{ll} 
u^L = u^R = 0, \\
\rho^{L,R} = \rho_{hs}^{L,R}, \\
p^{L,R} = p_{hs}^{L,R},
\end{array}\right.
\end{equation}
with $\rho_{hs}$ and $p_{hs}$ given hydrostatic states.
If the function $\bar{\rho}$ is defined by 
\begin{equation}
\label{eq:wb_rhobar_arbitrary}
\bar{\rho} (W^L,W^R) = \frac{1}{2} ( \rho^L+\rho^R )
\end{equation}
and the difference of the gravitational potential in the intermediate states is approximated by 
\begin{equation}
\label{eq:wb_ZR-ZL_approximation}
Z^R-Z^L \approx - \frac{p_{hs}^R-p_{hs}^L}{\frac{1}{2}(\rho_{hs}^L+\rho_{hs}^R)},
\end{equation}
then the approximate Riemann solver $W_\mathcal{R}$ preserves the steady state.
\end{lemma}

\begin{proof}
As can be seen in the proof of lemma \ref{lem:wb_3_families_hydrostatic_equilibria} it is sufficient to show that the conditions (\ref{eq:wb_steady_state_condition1}) and (\ref{eq:wb_steady_state_condition2}) are fulfilled so that lemma \ref{lem:wb_discrete_steady_states} can be applied.
In order to do so we plug the states from (\ref{eq:wb_arbitrary_hydrostatic_equilibrium}) and the approximation (\ref{eq:wb_ZR-ZL_approximation}) into (\ref{eq:wb_steady_state_condition2}) and use definition (\ref{eq:wb_rhobar_arbitrary}) for $\bar{\rho}$. This results in 
\begin{equation}
p_{hs}^R-p_{hs}^L-\frac{1}{2}(\rho_{hs}^L+\rho_{hs}^R)  \frac{p_{hs}^R-p_{hs}^L}{\frac{1}{2}(\rho_{hs}^L+\rho_{hs}^R)} = 0.
\end{equation}
\hfill
\end{proof}\\

Now that it has been shown that the approximate Riemann solver satisfies the well-balanced property, it remains to show that the entire scheme does so as well.

\setcounter{theorem}{16}
\begin{theorem}
\label{theo:wb}
Let us consider an initial data $\omega_i^0,\omega_{i+1}^0$ that satisfies
\begin{equation}
\begin{split}
u_i^0=u_{i+1}^0=0, \\
\frac{1}{\Delta x} (p_{i+1}^0 - p_i^0) + \bar{\rho}(W_i^0,W_{i+1}^0) \frac{\Phi_{i+1}-\Phi_i}{\Delta x} = 0.
\end{split}
\label{eq:wb_well-balanced_condition_in_the_scheme}
\end{equation}
Then the updated state $\omega^{n+1}$ stays at rest, and thus satisfies $\omega_i^{n+1} = \omega_i^n$ for
all $i\in \mathbb{Z}$.
\end{theorem}

\begin{proof}
Since both conditions (\ref{eq:wb_steady_state_condition1}) and (\ref{eq:wb_steady_state_condition2}) of lemma \ref{lem:wb_discrete_steady_states} are fulfilled, the approximate Riemann solver stays at rest. The updated state $\omega_i^1$ at time $t=\Delta t$ is in essence the sequence of approximate Riemann solvers. Since the approximate Riemann solver is at rest, it directly follows $\omega_i^1 = \omega_i^0$ for all $i\in \mathbb{Z}$.
\hfill
\end{proof}

\section{Extension to 2D}
\label{sec:2d_extension}
For two spatial dimensions the Euler equations \eqref{sys:Euler} can be written in the form
\begin{equation}
\omega_t + F(\omega)_x + G(\omega)_y = \mathcal{S}(\omega,\Phi).
\end{equation}
On a regular cartesian grid, we extend the numerical scheme described in Sect. \ref{sec:relaxation_scheme} to two dimensions by applying an unsplit finite volume method \cite{36}, in which the contributions of both directions are used in only one step to update the numerical solution by the formula
\begin{equation}
\label{eq:2d_relaxatation_scheme}
\begin{split}
\omega_{i,j}^{n+1} = \omega_{i,j}^n &- \frac{\Delta t}{\Delta x} \left( F_{i+1/2,j}^n - F_{i-1/2,j}^n \right) - \frac{\Delta t}{\Delta y} \left( G_{i,j+1/2}^n - G_{i,j-1/2}^n \right)\\ 
&+ \frac{\Delta t}{2} 
\left( S_{i-1/2,j}^{+,n} \frac{\Phi_{i,j}^n-\Phi_{i-1,j}^n}{\Delta x} + S_{i+1/2,j}^{-,n} \frac{\Phi_{i+1,j}^n-\Phi_{i,j}^n}{\Delta x} \right) \\
&+ \frac{\Delta t}{2} 
\left( S_{i,j-1/2}^{+,n} \frac{\Phi_{i,j}^n-\Phi_{i,j-1}^n}{\Delta x} + S_{i,j+1/2}^{-,n} \frac{\Phi_{i,j+1}^n-\Phi_{i,j}^n}{\Delta y} \right).
\end{split}
\end{equation}
The definitions of the numerical fluxes and source terms are straightforward extension of the ones in Sect. \ref{sec:relaxation_scheme}. The numerical fluxes continue to use the one-dimensional approximate Riemann solver, so that it is applied separately in $x$- and $y$-direction. This Riemann solver corresponds to that defined in \eqref{eq:approximate_riemann_solver}, in which addionally the intermediate states for the transversal velocity are set by the left and right values at the interface, respectively, since this velocity is a Riemann invariant for the outer waves $\sigma^-$ and $\sigma^+$.\\
Since the two-dimensional method is still based on the one-dimensional Riemann solver, the properties proven in Sect. \ref{sec:properties} also apply to this method. From this follows the entropy inequality, the absence of checkerboard modes, positivity and the asymptotic conservation property. In addition, the well-balanced property is also preserved, since the approximate Riemann solvers is at rest for initial data in hydrostatic equilibrium in both spatial directions and thus in both momentum equations the pressure gradient cancels out with the source term.


\section{Second order extension}
\label{sec:2nd_order}

In this section we give a possible extension of the proposed scheme to second order in space. 
We use a linear reconstruction in the primitive variables $\omega^p=(\rho,\boldsymbol{u},p)$. In order to obtain the values $\omega_{i-1/2}^R$ and $\omega_{i+1/2}^L$, which serve as initial data for the Riemann problems at the interface, we evaluate the function
\begin{equation}
\label{eq:2nd_order_linear_reconstruction}
\omega^p(x) = \omega_i^p + \sigma (x-x_i)
\end{equation}
in each cell $\mathcal{C}_i$ at its boundaries $x_{i-1/2}$ and $x_{i+1/2}$. The slope $\sigma$ depends on the neighbouring cells and is computed for each primitive variable separately.
In order to ensure that the reconstructed values for the density and internal energy remain positive, which is essential for the positivity property given by the lemmata \ref{lem:positivity} and \ref{lem:positivity_density}, we use a limiting procedure introduced in \cite{29} that builds on the work by Berthon in \cite{30}. Then the slopes are defined by
\begin{equation}
\label{eq:2nd_order_slopes_positivity}
\begin{split}
\sigma^{\rho} &= \rho_i \max \left ( -1, \min \left ( 1, \frac{\bar{\sigma}^\rho}{\rho_i} \right ) \right ), \\
\boldsymbol{\sigma^{u}} &= \kappa \boldsymbol{\bar{\sigma}^{u}}, \\
\sigma^{p} &= p_i \max \left ( -1, \min \left ( 1, \frac{\bar{\sigma}^p}{p_i} \right ) \right ),
\end{split}
\end{equation}
with
\begin{equation}
\label{eq:2nd_order_minmod}
\bar{\sigma} = \text{minmod} \left( \frac{\omega_{i}^p-\omega_{i-1}^p}{\Delta x}, \frac{\omega_{i+1}^p-\omega_i^p}{\Delta x}  \right)
\end{equation}
and
\begin{equation}
\label{eq:2nd_order_omega}
\begin{split}
\kappa &= \min (1, \bar{\kappa} ), \\
\bar{\kappa} &= 
\left\{\begin{array}{ll} 
\frac{-\sigma^\rho (\boldsymbol{u}_i \cdot \boldsymbol{\bar{\sigma}^{u}} ) + \sqrt{ (\sigma^\rho)^2 (\boldsymbol{u}_i \cdot \boldsymbol{\bar{\sigma}^{u}} )^2+
\Vert \boldsymbol{\bar{\sigma}^{u}} \Vert^2 \frac{\rho_i p_i}{\gamma-1}}}{ \rho_i 
\Vert \boldsymbol{\bar{\sigma}^{u}} \Vert^2}, &\text{if } \boldsymbol{\bar{\sigma}^{u}} \neq 0,  \\
1, &\text{if } \boldsymbol{\bar{\sigma}^{u}}=0.
\end{array}\right.
\end{split}
\end{equation}
Additionally, we also want to preserve the well-balanced property for the second-order scheme. To achieve this, we adjust the pressure slope by using a hydrostatic reconstruction \cite{28,29,8}. Instead of directly using the pressure values of the neighbouring cells, one first applies the transformations
\begin{equation}
\label{eq:2nd_order_transformation_q}
\begin{split}
q_{i-1}  &= p_{i-1} - \bar{\rho}(W_{i-1},W_i) (\Phi_i-\Phi_{i-1}), \\
q_{i+1} &= p_{i+1} + \bar{\rho}(W_i,W_{i+1}) (\Phi_{i+1}-\Phi_i),
\end{split}
\end{equation}
and then computes the slope for the pressure by 
\begin{equation}
\label{eq:2nd_order_slope_wb}
\bar{\sigma}^{p} = \text{minmod}\left( \frac{p_i-q_{i-1}}{\Delta x}, \frac{q_{i+1}-p_i}{\Delta x} \right).
\end{equation}
In the case of hydrostatic equilibrium, the slope becomes zero and the interface values for the pressure thus reduce to the cell averages.
The approximate Riemann solver then stays at rest due to lemma \ref{lem:wb_discrete_steady_states} and all results of the former section about well-balancing remain valid for the second order scheme.

The second order scheme remains asymptotic preserving since the differences in (\ref{eq:ap_numerical_derivative}) are due to the linear reconstruction of order $\mathcal{O}(\Delta x^2)$. In the following, we illustrate this for the velocity $u$ using backward slopes for $\sigma^{u}$:
\begin{align*}
u^L-u^R 
&= u_i+\sigma_i \frac{\Delta x}{2} - \left(u_{i+1}-\sigma_{i+1} \frac{\Delta x}{2} \right) \\
&= u_i + \frac{u_i-u_{i-1}}{\Delta x} \frac{\Delta x}{2} - u_{i+1} + \frac{u_{i+1}-u_i}{\Delta x} \frac{\Delta x}{2} \\
&= -\frac{1}{2} u_{i+1} + u_i - \frac{1}{2} u_{i-1} \\
&=- \frac{1}{2} (\Delta x)^2 \partial_x u.
\end{align*}
Thanks to these second order approximations of the derivatives, all first order terms in the consistency error can be replaced by second order terms so that the scheme becomes second order uniformly in $M$. The steps of proof of the second part of theorem \ref{theo:asymptotic_preserving} work analogously as for the first order case, if we assume for the last step the new condition $M<\sqrt{k \Delta x^2}$.


\section{Numerical results}
\label{sec:results}

In this section we numerically investigate the theoretical properties of the relaxation scheme presented in the previous sections. The approximate Riemann solver in the scheme is equipped with the intermediate states defined in (\ref{eq:IS_vs})-(\ref{eq:IS_abZ}) and the relaxation speeds (\ref{eq:ap_relaxation_speeds_with_fix}) with $\beta=1.1$.
Various definitions are used for the $\bar{\rho}$-function. Definition (\ref{eq:wb_rhobar_isothermal}) is used by default. If a different choice is made, this is indicated in the respective test. 
The second order spatial scheme is combined with a third order Runge Kutta method \cite{33} for time integration. 
For all test set-ups we assume an ideal gas law $p=(\gamma -1) \rho e$.
The computations are performed on a regular cartesian grid.

\subsection{Accuracy}
\label{sec:accuracy}
In a first numerical test, which is suggested by \cite{27}, we investigate the experimental order of convergence of the relaxation scheme presented. For the Euler equations (\ref{sys:Euler}) on the domain $[0,1]^2$ with a linear gravitational potential $\Phi(x,y)=x+y$, one possible exact solution is defined by
\begin{equation}
\label{eq:accuracy_exact_solution}
\begin{split}
\rho(x,y,t) =& 1 + 0.2 \sin \left( \pi \left(x + y - t(u_{1_0}+u_{2_0}) \right) \right), \\
\boldsymbol{u}(x,y,t) =& \left( u_{1_0},u_{2_0} \right), \\
p(x,y,t) =& 4.5 + (u_{1_0}+u_{2_0})t- (x+y) 
+ 0.2 \cos \left( (\pi \left(x + y - (u_{1_0}+u_{2_0})t \right) \right) / \pi,
\end{split}
\end{equation}
with $u_{1_0}=u_{2_0}=20$ and $p_0=4.5$. The exact solution is also used for the boundary conditions.  The adiabatic coefficient is set to $\gamma=5/3$. We compare the numerical and exact solutions computed on a $N\times N$ grid at final time $T=0.01$. The resulting $L^1$ errors and experimental orders of convergence can be found in Table \ref{tab:accuracy}. As expected, we obtain orders of convergence of nearly two. Without the use of limiters full second order is reached.

\begin{table}[h]
\scriptsize
\centering
 \begin{tabular}{c|cccccccc}
\specialrule{.2em}{.1em}{.1em}
$N$ & $L^1(\rho)$ & $EOC(\rho)$ & $L^1(\rho u_1)$ & $EOC(\rho u_1)$ & $L^1(\rho u_2)$ & $EOC(\rho u_2)$ & $L^1(E)$ & $EOC(E)$  \\
\specialrule{.05em}{.1em}{.1em} 
32   & 7.26E-04 &    -    & 1.45E-02 &    -    & 1.45E-02 &    -    & 2.90E-01 &     -   \\
64   & 1.97E-04 & 1.88 & 3.93E-03 & 1.88 & 3.93E-03 & 1.88 & 7.87E-02 & 1.88  \\ 
128 & 5.22E-05 & 1.92 & 1.04E-03 & 1.92 & 1.04E-03 & 1.92 & 2.08E-02 & 1.92  \\
256 & 1.37E-05 & 1.92 & 2.73E-04 & 1.93 & 2.73E-04 & 1.93 & 5.47E-03 & 1.93  \\
512 & 3.60E-06 & 1.94 & 7.10E-05 & 1.94 & 7.10E-05 & 1.94 & 1.42E-03 & 1.95   \\
\specialrule{.2em}{.1em}{.1em}
\end{tabular}
\caption{$L^1$ errors and experimental orders of convergence}
\label{tab:accuracy}
\end{table}

\subsection{Strong rarefaction test}
\label{sec:rarefaction}
\begin{figure}[t]
\begin{center}
\includegraphics[scale=0.6]{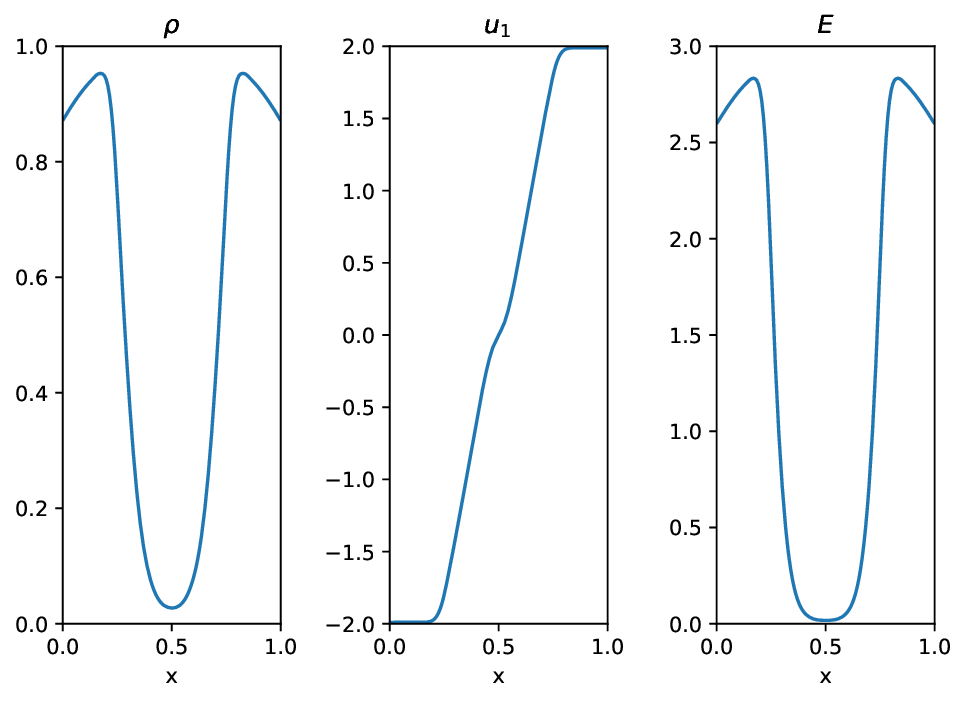}
\end{center}
 \caption{Numerical solution for density, velocity and total energy at final time $T=0.1$}
  \label{fig:rarefaction}
\end{figure}

In this section we want to numerically verify the theoretical results of Sect. \ref{sec:positivity}, i.e. the positivity of density and internal energy.
One suitable test for which density and pressure become very small is the 1-2-0-3 strong rarefaction test \cite{8}. In this test set-up, two rarefaction waves are launched in x-direction on top of an isothermal atmosphere. Therefore, on the domain $[0,1]^2$ the density $\rho$ and pressure $p$ are initially defined by (\ref{eq:wb_isothermal_equilibrium}) with the constants $C=-0.01$ and $K=\gamma-1$, an adiabatic coefficitent $\gamma=1.4$ and a quadratic gravitational potential 
$\Phi(x,y) = \frac{1}{2} [ (x-0.5)^2+(y-0.5)^2]$.
The initial velocities are set to
\begin{equation}
u_1 = 
\left\{\begin{array}{rl} -2, & x < 0.5, \\
							    2, & x \geq 0.5,
										\end{array}\right.
\quad \text{and} \quad
u_2=0.
\end{equation}
One slice along the $x$-axis of the numerical solution at final time $T=0.1$ computed on a $128\times 128$ grid by our relaxation scheme is presented in Fig. \ref{fig:rarefaction}. Although the values for density and total pressure become very small during the simulation, they always remain positive. This outcome underlines the theoretical results stated in lemmata \ref{lem:positivity} and \ref{lem:positivity_density}.

\subsection{Isothermal atmosphere}
\label{sec:isothermal_atmosphere}
The following set-up is taken from \cite{5}.
The aim of this experiment is to illustrate the exact preservation of an isothermal equilibrium. We consider the gravitational potential
\begin{equation}
\Phi(x,y) = x+y.
\end{equation}
The initial conditions on the domain $[0,1]^2$ are given by
\begin{equation}
\begin{split}
\rho(x,y,0) &= \rho_0 \exp(-\rho_0 g (x+y)/p_0),  \\
\boldsymbol{u}(x,y,0) &= 0, \\
 p(x,y,0) &= p_0 \exp(-\rho_0 g (x+y)/p_0),
\end{split}
\end{equation}
with the parameters $\rho_0=1.21$, $p_0=1$ and $g=1$. In this test we set  $\gamma=1.4$.
The solution should be preserved up to any final time. Here we choose $T=1.0$. 
Since the solution is in hydrostatic equilibrium, the choice of the $\bar{\rho}$-average plays an important role. As this is an isothermal equilibrium, we use for $\bar{\rho}$ the definition (\ref{eq:wb_rhobar_isothermal}). The $L^1$ error between the approximated solution and the exact solution is given in Table \ref{tab:isothermal_atmosphere} and is in the order of magnitude of the machine accuracy.

\begin{table}[h]
\scriptsize
\centering
 \begin{tabular}{c|cccccccc}
\specialrule{.2em}{.1em}{.1em}
$N$ & $L^1(\rho)$ & $L^1(\rho u_1)$ & $L^1(\rho u_2)$ & $L^1(E)$  \\
\specialrule{.05em}{.1em}{.1em} 
32   & 8.95E-17 & 5.21E-16 & 5.21E-16 & 4.18E-16 \\
64   & 1.73E-16 & 1.62E-16 & 1.62E-16 & 7.24E-16 \\
128 & 3.40E-16 & 3.47E-16 & 3.47E-16 & 1.63E-15 \\
256 & 6.30E-16 & 6.89E-16 & 6.89E-16 & 3.46E-15 \\
512 & 1.22E-15 & 1.54E-15 & 1.54E-15 & 7.43E-15 \\
\specialrule{.2em}{.1em}{.1em}
\end{tabular}
\caption{$L^1$ errors for an isothermal atmosphere}
\label{tab:isothermal_atmosphere} 
\end{table}

\subsection{General steady state}
\label{sec:general_steady_state}
In practice, steady states that do not belong to the class of polytropic equilibria can also occur. In order to investigate the behaviour of the well-balancing mechanism for these cases, we now apply the scheme to a general steady state. We take the initial conditions from the set-up in Sect. \ref{sec:accuracy} and set the initial velocities $u_{1_0}$ and $u_{2_0}$ to zero. Then it is easy to check that the initial data is in hydrostatic equilibrium.

In a first step, we use the $\bar{\rho}$-average tuned to isothermal equilibria (\ref{eq:wb_rhobar_isothermal}) and compute the solution at final time $T=1$. As expected, the $L^1$ error shown in Table \ref{tab:general_steady_state_isothermal_rhobar} is now no longer in the order of magnitude of the machine accuracy, but the hydrostatic equilibrium is still preserved up to second order. This result remains true even if we use a constant reconstruction and consequently a first order scheme. As the convergence rates in Table \ref{tab:general_steady_state_isothermal_rhobar_1st_order} show, the hydrostatic equilibrium is maintained up to second order despite the constant reconstruction. Mathematically, this can be explained by the fact that equation (\ref{eq:wb_well-balanced_condition_in_the_scheme}) is satisfied up to second order.

\begin{table}[h]
\scriptsize
\centering
 \begin{tabular}{c|cccccccc}
\specialrule{.2em}{.1em}{.1em}
$N$ & $L^1(\rho)$ & $EOC(\rho)$ & $L^1(\rho u_1)$ & $EOC(\rho u_1)$ & $L^1(\rho u_2)$ & $EOC(\rho u_2)$ & $L^1(E)$ & $EOC(E)$  \\
\specialrule{.05em}{.1em}{.1em} 
32   & 9.43E-06 &    -    & 1.36E-05 &    -    & 1.36E-05 &    -    & 5.08E-05 &     -   \\
64   & 2.35E-06 & 2.01 & 3.43E-06 & 1.99 & 3.43E-06 & 1.99 & 1.26E-05 & 2.01  \\ 
128 & 5.88E-07 & 2.00 & 8.60E-07 & 2.00 & 8.60E-07 & 2.00 & 3.14E-06 & 2.01  \\
256 & 1.47E-07 & 2.00 & 2.16E-07 & 1.99 & 2.16E-07 & 1.99 & 7.85E-07 & 2.00  \\
512 & 3.69E-08 & 2.00 & 5.42E-08 & 2.00 & 5.42E-08 & 2.00 & 1.97E-07 & 2.00  \\
\specialrule{.2em}{.1em}{.1em}
\end{tabular}
\caption{$L^1$ errors and experimental orders of convergence of the second order scheme for a general steady state using the $\bar{\rho}$-average (\ref{eq:wb_rhobar_isothermal})}
\label{tab:general_steady_state_isothermal_rhobar}
\end{table}

\begin{table}[h]
\scriptsize
\centering
 \begin{tabular}{c|cccccccc}
\specialrule{.2em}{.1em}{.1em}
$N$ & $L^1(\rho)$ & $EOC(\rho)$ & $L^1(\rho u_1)$ & $EOC(\rho u_1)$ & $L^1(\rho u_2)$ & $EOC(\rho u_2)$ & $L^1(E)$ & $EOC(E)$  \\
\specialrule{.05em}{.1em}{.1em} 
32   & 9.74E-06 &    -    & 1.40E-05 &    -    & 1.40E-05 &    -    & 5.15E-05 &     -   \\
64   & 2.39E-06 & 2.03 & 3.48E-06 & 2.01 & 3.48E-06 & 2.01 & 1.27E-05 & 2.02  \\ 
128 & 5.93E-07 & 2.01 & 8.67E-07 & 2.01 & 8.67E-07 & 2.01 & 3.15E-06 & 2.01  \\
256 & 1.48E-07 & 2.00 & 2.17E-07 & 2.00 & 2.17E-07 & 2.00 & 7.86E-07 & 2.00  \\
512 & 3.70E-08 & 2.00 & 5.43E-08 & 2.00 & 5.43E-08 & 2.00 & 1.97E-07 & 2.00   \\
\specialrule{.2em}{.1em}{.1em}
\end{tabular}
\caption{$L^1$ errors and experimental orders of convergence of the first order scheme for a general steady state using the $\bar{\rho}$-average (\ref{eq:wb_rhobar_isothermal})}
\label{tab:general_steady_state_isothermal_rhobar_1st_order}
\end{table}

\noindent
Let us now assume that we know the hydrostatic equilibrium a priori and it is given as discrete data for the density and pressure. In this case, the approach described in lemma \ref{lem:wb_arbitrary_hydrostatic_equilibrium} should be able to maintain this particular hydrostatic equilibrium up to machine precision. In order to check this, we set the values $\rho_{hs}$ and $p_{hs}$ equal to the initial values for density respective pressure. The $L^1$ errors in Table \ref{tab:general_steady_state_arbitrary_rhobar} show that the hydrostatic equilibrium is maintained up to machine precision.

\begin{table}[h]
\scriptsize
\centering
 \begin{tabular}{c|cccccccc}
\specialrule{.2em}{.1em}{.1em}
$N$ & $L^1(\rho)$ & $L^1(\rho u_1)$ & $L^1(\rho u_2)$ & $L^1(E)$  \\
\specialrule{.05em}{.1em}{.1em} 
32   & 6.54E-17 & 9.10E-16 & 9.10E-16 & 1.33E-15 \\
64   & 1.85E-16 & 1.95E-15 & 1.95E-15 & 4.78E-15 \\
128 & 2.98E-16 & 4.78E-15 & 4.78E-15 & 8.97E-15 \\
256 & 6.25E-16 & 8.32E-16 & 8.32E-16 & 2.04E-14 \\
512 & 1.25E-15 & 1.83E-14 & 1.83E-14 & 4.24E-14 \\
\specialrule{.2em}{.1em}{.1em}
\end{tabular}
\caption{$L^1$ errors for a general steady state using the approach for a-priori known hydrostatic equilibria from lemma \ref{lem:wb_arbitrary_hydrostatic_equilibrium}}
\label{tab:general_steady_state_arbitrary_rhobar}
\end{table}

\subsection{Perturbation of an isothermal atmosphere}
\label{sec:perturbation}
One main advantage of well-balanced schemes is their ability to resolve small perturbations on the hydrostatic equilibrium even on coarse grids. It is precisely this effect that we are investigating with the following test.
For this purpose, we take the initial values from Sect. \ref{sec:isothermal_atmosphere}, which are in hydrostatic equilibrium, and add a perturbation on the pressure
\begin{equation}
p(x,y,0) = p_0 \exp \left(-\rho_0 g (x+y)/p_0 \right) 
+ \eta \exp \left(-100 \rho_0 g ((x-0.3)^2+(y-0.3)^2)/p_0 \right).
\end{equation}
The strength of the perturbation is controlled by the parameter $\eta$. The numerical solutions are computed on a $64\times 64$ mesh up to a final time $t=0.15$. 
In order to investigate the well-balancing effect, we compare the results of our well-balanced scheme with a non-well-balanced scheme. The non-well-balanced scheme uses a Rusanov flux in combination with a linear reconstruction limited by the minmod limiter.

\begin{figure}[t]
\centering
\includegraphics[scale=0.7]{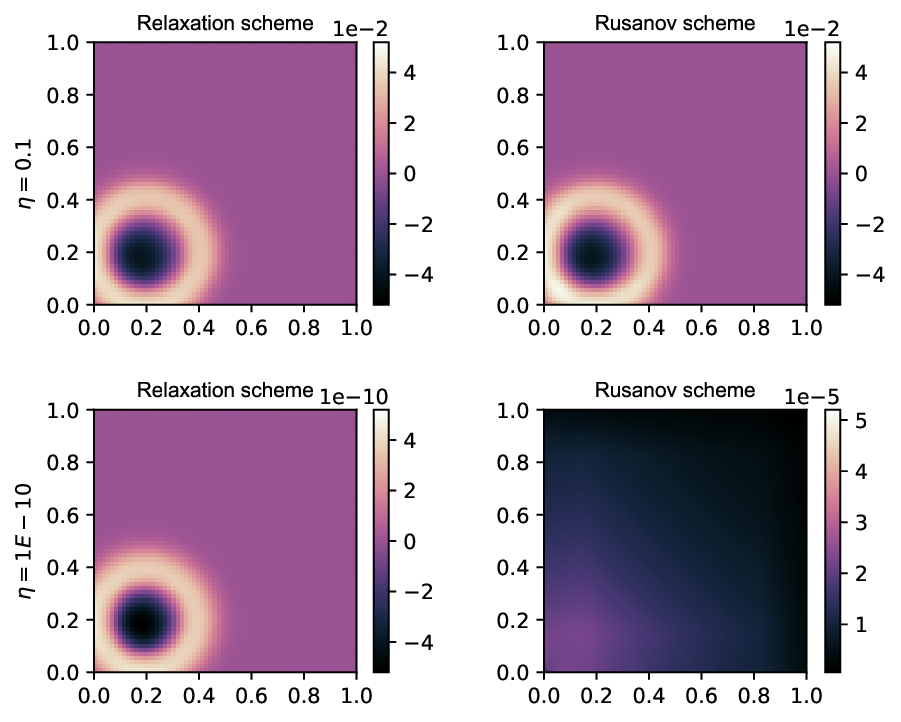}
\caption{Pressure perturbation of an isothermal atmosphere at $t=0.15$}
 \label{fig:perturbation}
\end{figure}

The numerical solutions of the two schemes for a large perturbation $\eta=0.1$ are illustrated in the two upper plots of Fig. \ref{fig:perturbation}. Looking at the two solutions, it can be said that they are visually very similar. Both methods are capable of resolving the perturbation well.
For a significantly smaller perturbation ($\eta=1$E$-10$), the situation is completely different. While our well-balanced relaxation scheme is still capable to resolve the perturbation (in fact one cannot see any difference in the resolution in comparison to the larger perturbation), the non-well-balanced scheme completely destroys the structure of the initial pressure pulse. This underlines the functionality of the well-balanced mechanism in the relaxation scheme and also demonstrates the importance of this property for problems close to hydrostatic equilibrium.

\subsection{Stationary vortex in a gravitational field}
\label{sec:low_mach_vortex}
In this section, we investigate the effect of the new scaling of the relaxation speeds $a$ and $b$ for problems with low Mach numbers. Therefore we compare the two-speed relaxation scheme using the speeds defined in (\ref{eq:ap_relaxation_speeds_with_fix}) with the one-speed relaxation scheme using the speeds (\ref{eq:ap_relaxation_speeds_without_fix}).

\begin{figure}[t]
\begin{center}
\includegraphics[scale=0.8]{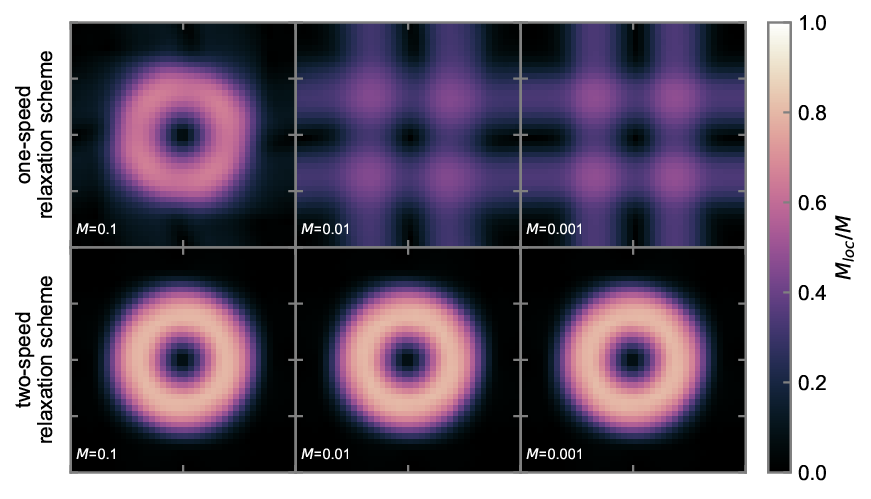}
\end{center}
  \caption{Numerical solutions for different maximum Mach numbers $M$ after one full turnover. The local Mach number relative to the respective $M$ is color coded}
 \label{fig:vortex}
\end{figure}

As a test, we use a version of the Gresho vortex modified for the Euler equations with a gravitational source term that was already given in \cite{8}.
The density in this set-up is defined by
\begin{equation}
\rho = \exp \left( -\frac{\Phi}{RT} \right).
\end{equation}
The rest of the initial data is given in radial coordinates $(r,\theta)$. The velocity field has the form
\begin{equation}
u_\theta(r) = \frac{1}{u_r}
\left\{\begin{array}{ll} 5r, & r \leq 0.2 \\
							  2-5r, & 0.2< r \leq 0.4 \\
							 0,  & r>0.4
										\end{array}\right.
\end{equation}
and the gravitational potential is defined by
\begin{equation}
\Phi (r) = 
\left\{\begin{array}{ll} 12r^2, & r \leq 0.2 \\
							  0.5-\ln (0.2)+ \ln (r), & 0.2< r \leq 0.4 \\
							 \ln(2)-0.5 \frac{r_c}{r_c-0.4}+2.5\frac{r_c}{r_c-0.4}r-1.25\frac{1}{r_c-0.4}r^2,  & 0.4<r\leq r_c \\
							 \ln(2)-0.5 \frac{r_c}{r_c-0.4}+1.25\frac{r_c^2}{r_c-0.4}, & r>r_c.
										\end{array}\right.
\end{equation}
The pressure $p$ is departed into a hydrostatic pressure $p_0$ and a pressure $p_2$ associated with the centrifugal forces and given by $p=p_0+M^2 p_2$, where $p_0=RT \rho$ and
\begin{equation}
p_2(r) = \frac{RT}{u_r^2}
\left\{\begin{array}{ll} p_{21}(r), & r \leq 0.2 \\
							  p_{21}(0.2)+p_{22}(r), & 0.2< r \leq 0.4 \\
							  p_{21}(0.2)+p_{22}(0.4),  & 0.4<r\leq r_c
										\end{array}\right.
\end{equation}
with
\begin{align*}
p_{21}(r) =& \left( 1 - \exp \left( -12.5 \frac{r^2}{RT} \right) \right), \\
\begin{split}
p_{22}(r) =& \frac{1}{(1-M^2) (1-0.5M^2)} \exp \left( \frac{-0.5+\ln(0.2))}{RT}  \right) \\
&\bigg ( r^{-\frac{1}{RT}} \left( M^4 (r (10-12.5r) - 2 ) - 4 + M^4 (r (12.5r - 20 ) + 6 ) RT \right) \\
&+ \exp \left( \frac{-\ln(0.2)}{RT} \right) \left( 4 - 2.5 M^4 RT + 0.5 M^4 \right) \bigg ).
\end{split}
\end{align*}
The reference values are given by $u_r=2\cdot 0.2 \cdot \pi$ and $RT=1/M^2$. We choose $\gamma=5/3$ for the adiabatic coefficient. The spatial domain is $D=[0,1]^2$ and has periodic boundary conditions. The computations are carried out on a $40 \times 40$ grid until a final time $T=1$, which corresponds to one turn of the vortex.
We solve this initial value problem for various maximum Mach numbers $M$ using the two different schemes. The solutions generated by the one-speed relaxation scheme are depicted in the top row of Fig. \ref{fig:vortex}, while the solutions computed by the two-speed relaxation scheme are shown in the bottom row. It becomes clear that for decreasing Mach numbers, the vortex in the upper row smears out very quickly and even loses its shape completely. The vortex produced by the two-speed scheme in the lower row, on the other hand, retains its shape regardless of the Mach number, so that no difference is visually discernible.

This outcome can be explained by the theoretical results from Sect. \ref{sec:AP}. While the diffusion for the one-speed scheme increases for decreasing Mach numbers, the use of two relaxation speeds, as shown in the proof of theorem \ref{theo:asymptotic_preserving}, results in a Mach number independent diffusion. Further evidence for this behaviour can be found in the analysis of the kinetic energy. Table \ref{tab:kinetic_energy} contains the percentages of kinetic energy compared to the initial value after one full turnover. The final amount of kinetic energy in the solutions of the one-speed scheme strongly decreases for decreasing Mach numbers. In contrast, the asymptotic preserving two-speed scheme is able to keep the loss almost constant at 14\%, regardless of the Mach number.

\begin{table}[h]
\scriptsize
\centering
 \begin{tabular}{c|ccc}
\specialrule{.2em}{.1em}{.1em}
scheme & $M=0.1$ & $M=0.01$ & $M=0.001$  \\
\specialrule{.05em}{.1em}{.1em} 
1-speed   & 62.74 & 47.49  & 50.08 \\
2-speed   & 86.03 & 86.00  & 85.99  \\
\specialrule{.2em}{.1em}{.1em}
\end{tabular}
\caption{Percentage of kinetic energy compared to the initial value after one full turnover ($T=1$) for the stationary vortex in a gravitational field}
\label{tab:kinetic_energy}
\end{table}


\section{Conclusion}
\label{sec:conclusion}
The proposed scheme extends the two-speed relaxation approach to the full Euler equations with a gravitational source term. In order to preserve steady states at rest, a well-balancing mechanism is installed in the approximate Riemann solver. The resulting scheme is provably asymptotic preserving and maintains all hydrostatic equilibria up to second order, certain families and a-priori known equilibria even up to machine precision. The approximate Riemann solver is positivity preserving, entropy satisfying and prevents the occurrence of checkerboard modes in the velocity and pressure variables. The properties of the method proven in theory are substantiated in numerical tests.
Further steps may be the development of an IMEX scheme based on the herein presented full time explicit scheme in order to overcome the severe time step restriction for problems with very low Mach numbers, and the extension to other PDE systems like the MHD equations.

\section*{Acknowledgements}
The authors thank Jonas Berberich for his comments on the well-balancing of a-priori known hydrostatic equilibria. We acknowledge the use of the \textit{Seven-League Hydro Code} (\url{https://slh-code.org}) for our numerical experiments. Claudius Birke acknowledges the support by the German Research Foundation (DFG) under the project no. KL 566/22-1. All three authors acknowledge the project  "Bayerisch-Franz\"osische Hochschulzentrum FK40\_2019'' which supported this work.




\end{document}